\documentclass[11pt,a4paper]{article}
\IfFileExists{lmodern.sty}{\usepackage[T1]{fontenc}\usepackage{lmodern}}{\usepackage[OT1]{fontenc}}
\usepackage[expansion=false]{microtype}
\usepackage{amsmath,amssymb,amsthm,mathtools}
\usepackage{geometry}
\usepackage{enumitem}
\usepackage{hyperref}
\usepackage{xcolor}
\hypersetup{colorlinks=true,linkcolor=blue!50!black,citecolor=blue!50!black,urlcolor=blue!50!black,pdftitle={Separable models of an NIP Banach space theory},pdfauthor={Karim Khanaki and Pedro Poitevin},pdfkeywords={continuous logic, Banach spaces, NIP, Hilbert renorming}}
\allowdisplaybreaks
\newtheorem{theorem}{Theorem}[section]
\newtheorem{lemma}[theorem]{Lemma}
\newtheorem{proposition}[theorem]{Proposition}
\newtheorem{corollary}[theorem]{Corollary}
\theoremstyle{definition}

\newcommand{\Hil}{\mathcal H}
\newcommand{\Y}{Y}
\newcommand{\Nrm}{\mathsf N}
\newcommand{\h}{\mathsf h}
\newcommand{\clsp}{\overline{\operatorname{span}}}
\newcommand{\supp}{\operatorname{supp}}
\newcommand{\dist}{\operatorname{dist}}

\newcommand{\aut}{\operatorname{Aut}}
\newcommand{\cU}{\mathcal U}
\newcommand{\ind}{\mathbf 1}
\title{\bfseries Separable models of an NIP Banach space theory}
\author{%
  \textbf{Karim Khanaki}\\[0.3em]
  {\small Department of Science}\\
  {\small Arak University of Technology}\\
  {\small Arak 38181-15881, Iran}\\[0.2em]
  {\footnotesize \href{mailto:khanaki@arakut.ac.ir}{\texttt{khanaki@arakut.ac.ir}}}
  \and
  \textbf{Pedro Poitevin}\\[0.3em]
  {\small Department of Mathematics}\\
  {\small Salem State University}\\
  {\small 352 Lafayette Street, Salem, MA 01970, USA}\\[0.2em]
  {\footnotesize \href{mailto:ppoitevin@salemstate.edu}{\texttt{ppoitevin@salemstate.edu}}}
}
\date{}
\begin{document}
\maketitle
\begin{abstract}
\noindent
Let $Y$ be the renorming of $\ell_2(\mathbb N)\oplus_2\ell_2(\mathbb N)$
with norm $(\h^2+q^2)^{1/2}$, where $\h$ is the Hilbert norm and
$q(u,v)=\sup_{i<j}|u_i+v_j|$. We prove that the unstable theory
$\operatorname{Th}(Y)$ is NIP in the pure Banach language. We recover
the index order from the Banach norm and classify the continuum many
separable models up to surjective linear isometry. The invariants are
the order type of a countable model of $\operatorname{Th}(\mathbb N,<)$
and the dimension of a separable Hilbert summand. Every such pair
occurs. We describe all elementary embeddings between these models.
The bounded norm types over a separable model are separable in the
uniform metric exactly when its index order contains no copy of
$(\mathbb Q,<)$.
\end{abstract}

\section{Introduction and preliminaries}
\subsection{Introduction}
The uniformly convex renorming $Y$ of Hilbert space was introduced in
\cite[Section~4.1 and Theorems~1.1 and~4.7]{PoitevinNormFormula2026v2},
where it was shown that $\operatorname{Th}(Y)$ is unstable and that
the norm formula has NIP in this theory.
We prove NIP for every formula in the pure Banach language and
classify the separable models of $\operatorname{Th}(Y)$ and the
elementary embeddings between them.

Write $\mathbb N=\{1,2,\dots\}$ and
$\Hil=\ell_2(\mathbb N)\oplus_2\ell_2(\mathbb N)$ over the reals.
For $z=(u,v)\in\Hil$, let $\h(z)$ be its Hilbert norm and define
\begin{equation}\label{eq:norm}
 q(z)=\sup_{i<j}|u_i+v_j|,
 \qquad \Nrm(z)=\bigl(\h(z)^2+q(z)^2\bigr)^{1/2}.
\end{equation}
$q$ is a seminorm because it is the supremum of absolute values of
linear functionals. We have $q(u,v)\le\|u\|_2+\|v\|_2\le\sqrt2\,\h(u,v)$,
so $\h\le\Nrm\le\sqrt3\,\h$. Put $Y=(\Hil,\Nrm)$. The same seminorm
occurs in the renorming $\h+\varepsilon q$ of $\ell_2$ for $\varepsilon>0$.
Iovino gives the formula and, following Garling, credits Bollob\'as with the example
\cite[1995 author version, p.~13, after Theorem~8.7]{IovinoStableBanachI}%
~\cite{Garling1982}.

Let $e_i$ be the $i$th unit vector of $\ell_2(\mathbb N)$.
Then $q(e_i,e_j)=2$ if $i<j$. If $i\ge j$, at most one nonzero
coefficient occurs in each sum and the pair $(i,i+1)$ gives a sum
of $1$. Since $\h(e_i,e_j)^2=2$, we have
\[
 \Nrm((e_i,e_j))=
 \begin{cases}
   \sqrt6,& i<j,\\
   \sqrt3,& i\ge j.
 \end{cases}
\]
The iterated limits are $\sqrt6$ with $j\to\infty$ first and $\sqrt3$
with $i\to\infty$ first, so $\Nrm(x+y)$ fails the double-limit condition
for stability \cite{KrivineMaurey} on $Y$. In particular,
$\operatorname{Th}(Y)$ is unstable.

A formula $\varphi(x;y)$ on fixed ball sorts has no independence
property (NIP) in a theory if, for every $r<s$, there is a bound on the
length of a list $a_1,\ldots,a_n$ in any model of the theory for which
varying $b$ realizes all $2^n$ patterns of values $\varphi(a_i;b)\le r$
or $\varphi(a_i;b)\ge s$. Theorem~\ref{thm:main} establishes this
bound for every formula in the pure Banach language of $Y$.

In \cite[Theorems~9.4 and~9.5]{LupiniNormFormulas2026},
Lupini proves NIP for a Hilbert space expanded by an index sort
$(\mathbb Q,<)$ and maps naming orthonormal families. He deduces
NIP for the complete Banach theory of an additive renorming defined
using ordered pairs of indices. His approximation of formula values
uses large coefficients, the order of their indices, and the full
Hilbert Gram matrix. In the discrete order used here, matching the
indices of large coefficients also requires comparing the number of
points before the first index and between successive indices, up to
a cutoff (Lemma~\ref{lem:formula-approximation}). For $n$ fixed
tuples, simultaneous approximation of the Hilbert inner products
with a parameter, together with the coefficient and gap comparisons,
gives a polynomial bound in $n$ on the number of patterns between
fixed thresholds $r<s$ (Proposition~\ref{prop:determined-subsets}).
We obtained an initial proof that $\operatorname{Th}(Y)$ is NIP before
learning of Lupini's preprint.

The order $(\mathbb N,<)$ used in the NIP proof is interpretable
without parameters in the pure Banach structure $Y$
(Theorem~\ref{thm:order-interpretation}). The interpretation uses a
relation on the pairs $i<j$ that identifies when two pairs share
a row or a column. The separable models of
$\operatorname{Th}(Y)$ are, up to surjective linear isometry,
precisely the spaces $Y_I^\circ\oplus_2 K$, where $I$ is a countable
model of $\operatorname{Th}(\mathbb N,<)$ and $K$ is a separable real
Hilbert space, possibly zero. Here $Y_I^\circ$ is the renorming of
$\ell_2(I)\oplus_2\ell_2(I\setminus\{\min I\})$ given by
\eqref{eq:norm} with the order on $I$. Two such models are
surjectively linearly isometric exactly when their interpreted
orders are isomorphic and their Hilbert summands have the same
dimension (Theorem~\ref{thm:separable-classification}).

Elementary embeddings between these models correspond bijectively
to triples consisting of an elementary embedding of the orders,
a common sign on the coordinate summand, and a linear isometric
embedding between the Hilbert summands
(Theorem~\ref{thm:embedding-normal-form}).

Every countable linear order $B$ occurs in the classification. The
corresponding index order consists of an initial copy of $\omega$
followed by copies of $\mathbb Z$ arranged in the order $B$.
There are exactly continuum
many separable models up to surjective linear isometry
(Corollary~\ref{cor:separable-spectrum}). The prime model is
$Y_{\mathbb N}^\circ$, whose image in a monster model is the
algebraic closure of the empty set in the Banach sorts
(Corollaries~\ref{cor:prime} and~\ref{cor:algebraic-closure}).

For a separable model $M$, the bounded norm types form the pointwise
closure on its closed unit ball of the functions $x\mapsto\|x+y\|_M$ with $y$
in that ball. The types are \emph{strongly separable} if this space
is separable in the uniform metric. This occurs exactly
when the interpreted order contains no copy of $(\mathbb Q,<)$
(Theorem~\ref{thm:strong-types}). The bounded norm types over $Y$ are
therefore separable in the uniform metric. A separable elementary
extension of $Y$ with index order $\omega+\sum_{r\in\mathbb Q}\mathbb Z$
has continuum many bounded norm types at pairwise distance bounded
below by a positive constant.

The theory fails quantifier elimination even after naming $\h$ and $q$
(Proposition~\ref{cor:no-QE}). An explicit quantifier-free Banach formula
witnesses the strict order property in the sense
of \cite[Definition~3.4]{Khanaki2026}
(Theorem~\ref{thm:SOP}).

\subsection{An auxiliary structure}
We regard $Y$ as a structure in the standard continuous Banach
language, with a sort for each closed ball of nonnegative rational
radius. In the underlying Hilbert space $\Hil$, write
\[
  E_i=(e_i,0),\quad F_i=(0,e_i).
\]
An auxiliary continuous structure $\mathcal M$ has the following sorts,
operations, and maps.
\begin{enumerate}[label=(\roman*),leftmargin=3em]
\item The Hilbert sort $H$ is realized in $\mathcal M$ by the real
Hilbert space $\Hil$, with its usual bounded ball sorts, inner
product, rational linear operations, and norm.
\item The discrete index sort $I=(\mathbb N,<)$ has metric $d_I(i,j)=\ind_{i\ne j}$.
\item The uniformly continuous maps $E,F:I\to H$ are given by
$E(i)=E_i$ and $F(i)=F_i$.
\end{enumerate}
Write $\chi_{<}(i,j)=1$ if $i<j$ and $\chi_{<}(i,j)=0$ otherwise.
On each bounded Hilbert ball, write $q$ for the formula
\begin{equation}\label{eq:qdef}
 q(z)=\sup_{i,j\in I}\chi_{<}(i,j)
      \big|\langle z,E(i)\rangle+\langle z,F(j)\rangle\big|.
\end{equation}
In $\mathcal M$ this formula agrees with the seminorm in \eqref{eq:norm}.
It also defines $q$ in elementary extensions of $\mathcal M$, where the
supremum is over the corresponding index sort. The norm $\Nrm$
is then defined from $\h$ and $q$ by \eqref{eq:norm}.

\subsection{Logical conventions and Banach-ball interpretation}
We use continuous first-order logic with bounded metric sorts
\cite{BBHU08}. We use continuous connectives on fixed bounded ranges.
The index sort has the $0$--$1$ metric and a two-valued order
predicate, so it can be viewed as a classical linear order in every
elementary extension.
The Banach language is countable. Quantifiers in Banach formulas
range over $\Nrm$-balls. We interpret each such ball inside the
Hilbert ball of the same radius.

\begin{lemma}[Uniform interpretation of Banach balls]\label{lem:interp}
For each rational $R\ge0$, the $\Nrm$-ball
$B_R^{\Nrm}=\{z:\Nrm(z)\le R\}$ is a definable subset of the Hilbert ball
$B_R^\h=\{z:\h(z)\le R\}$. With the metric $\Nrm(z-w)$ and the
restrictions of the Hilbert operations, these balls interpret the
Banach structure uniformly in models of $\operatorname{Th}(\mathcal M)$.
The interpretation
commutes with elementary embeddings and ultraproducts.
\end{lemma}
\begin{proof}
On $B_R^\h$ put $p_R(z)=(\Nrm(z)-R)_+$. If $\Nrm(z)>R$, the vector
$w=Rz/\Nrm(z)$ lies in $B_R^{\Nrm}$, and $\h(z-w)\le p_R(z)$ since
$\h(z)\le\Nrm(z)$. If $\Nrm(z)\le R$, take $w=z$. In either case,
$\dist_\h(z,B_R^{\Nrm})\le p_R(z)$.

For a $[0,1]$-valued ambient definable predicate $\psi(z,\bar a)$
on fixed bounded sorts, let $\omega$ be a modulus of continuity
in $z$. Then
\[
 \inf_{z\in B_R^\h}\min\{1,\psi(z,\bar a)+k p_R(z)\}
 \ \xrightarrow{k\to\infty}\ \inf_{z\in B_R^{\Nrm}}\psi(z,\bar a).
\]
The left side is less than or equal to the right side, since
$p_R$ vanishes on $B_R^{\Nrm}$. If $p_R(z)\ge1/k$, the expression
inside the infimum is $1$. Otherwise $z$ is within Hilbert distance
$1/k$ of its radial projection into $B_R^{\Nrm}$. Hence the error is at
most $\omega(1/k)$, uniformly in the models and the parameters.
Taking $\psi(z,a)$ to be a rescaled Hilbert distance to $a$ gives
definability of the distance to $B_R^{\Nrm}$ on $B_R^\h$. Suprema follow by
complementation. Translating atomic formulas by the Hilbert operations
and $\Nrm$ and applying the quantifier translations inductively gives
the uniform interpretation and preservation under elementary embeddings.

By $\h\le\Nrm\le\sqrt3\,\h$, both norms give the same bounded
families and the same families whose norms tend to zero along
$\cU$. Thus the two ultraproducts have the same vector classes
\cite{Heinrich}.
A class in the interpreted $\Nrm$-ball has representatives
$z_s\in B_R^\h$ with $\lim_{\cU}p_R(z_s)=0$.
Projecting each $z_s$ radially into $B_R^{\Nrm}$ moves it by at most
$p_R(z_s)$ in Hilbert norm, so the projected family represents the
same class. Conversely, a family with each $z_s\in B_R^{\Nrm}$ has
$\lim_{\cU}\Nrm(z_s)\le R$. The norm and vector operations agree
on these representatives, so the interpreted structure is the
Banach ultraproduct.
\end{proof}

\section{NIP of the complete theory}\label{sec:NIP}
\subsection{Order types and Hilbert tuples}
Let $\mathcal M^*\succ \mathcal M$ be an elementary extension with sorts $H^*$ and
$I^*$. By elementarity, the family
$\{E_\alpha,F_\alpha:\alpha\in I^*\}$ is orthonormal, and
\begin{equation}\label{eq:qstar}
 q^*(z)=\sup_{\alpha<\beta\ \text{in }I^*}
   |\langle z,E_\alpha\rangle+\langle z,F_\beta\rangle|.
\end{equation}
Set
\[
 A^*=\clsp\{E_\alpha,F_\alpha:\alpha\in I^*\},\qquad
 D^*=(A^*)^\perp.
\]
The orthogonal decomposition $H^*=A^*\oplus D^*$ holds for
any cardinality of $I^*$.
Since all indexed vectors lie in $A^*$, formula \eqref{eq:qstar}
gives $q^*(z)=q^*(P_{A^*}z)$.
For a finite Hilbert tuple $\bar a=(a_1,\dots,a_m)$, define its
\emph{indexed support} by
\[
 J(\bar a)=\bigl\{\alpha\in I^*: \exists r\le m,
   \langle a_r,E_\alpha\rangle\ne0
   \text{ or } \langle a_r,F_\alpha\rangle\ne0\bigr\}.
\]
For every $t>0$ the set of indices with a coefficient of magnitude
at least $t$ is finite by Bessel's inequality
\cite[Corollaries~4.9--4.10]{Conway1990}
applied to the orthonormal family $\{E_\alpha,F_\alpha:\alpha\in I^*\}$. Taking $t=1/n$ shows that
$J(\bar a)$ is countable, even if $I^*$ is uncountable.

\subsubsection{Types of finite tuples in discrete linear orders}
In a linear order, the \emph{initial gap} of an increasing tuple
$(a_1<\cdots<a_r)$ is the number of points preceding $a_1$. Its
successive gaps are the numbers of points strictly between $a_i$
and $a_{i+1}$. We allow infinite gaps. For an integer $L\ge0$,
two gaps have equal $L$-truncations if both are at least $L$ or
if both have the same length less than $L$.

\begin{lemma}\label{lem:ord}
Let $I$ and $J$ be linear orders with a least element, no greatest
element, an immediate successor for each point, and an immediate
predecessor for each nonleast point. This includes any two models
of $\operatorname{Th}(\mathbb N,<)$. Let $\bar a$ and $\bar b$
be increasing nonempty finite tuples of the same length in $I$
and $J$, respectively.
\begin{enumerate}[label=(\alph*),ref=\alph*,leftmargin=3em]
\item\label{item:ord-rank} For each $d\ge0$, if the initial and
successive gaps agree after truncation at $2^d-1$, the tuples agree
on all pure order formulas of quantifier rank at most $d$.
\item\label{item:ord-types} The tuples have the same complete order
type over $\varnothing$ if and only if the corresponding gaps have
equal $L$-truncations for every positive integer $L$.
\item\label{item:ord-equivalence} The least elements of $I$ and $J$
have the same complete order type, and $I\equiv J$.
\end{enumerate}
\end{lemma}
\begin{proof}
For (\ref{item:ord-rank}), play the Ehrenfeucht--Fra\"iss\'e game
\cite[Sections~3.2--3.3]{Hodges1993} with the two tuples already paired.
For each $d\ge0$, put
\begin{equation}\label{eq:threshold}
 T_d=2^d-1.
\end{equation}
Then $T_{d+1}=2T_d+1$. With $d$ moves remaining, we require
corresponding gaps to have equal lengths less than $T_d$ or
to both have length at least $T_d$.
With $d+1$ moves remaining, if the chosen point
lies in a gap of length less than
$T_{d+1}$, choose the point at the same position from the left.
The two new gaps have the same lengths as their counterparts.
For a gap
of length at least $T_{d+1}$, let $\ell$ and $r$ be the lengths of
the gaps to the left and right of the chosen point. At most one
of $\ell$ and $r$ is less than $T_d$, since otherwise
$\ell+r+1\le 2T_d-1$.
\begin{itemize}
\item If $\ell<T_d$, Duplicator chooses the point with exactly
$\ell$ points to its left in the corresponding gap. Both gaps to
the right of the paired points have at least $T_d$ points, since
each original gap has at least $2T_d+1$ points.
\item If $r<T_d$, Duplicator chooses the point with exactly
$r$ points to its right. The same count gives at least $T_d$
points to the left of each paired point.
\item If $\ell,r\ge T_d$, Duplicator chooses the $(T_d+1)$st
point from the left in the corresponding gap. Its length is at least
$2T_d+1$, so the response leaves at least $T_d$ points on both sides.
\end{itemize}

The same split works in the initial gap, using the least element
as the first available position. If the chosen point lies in the final
ray, count the points strictly between it and the greatest point already
matched in that order. Choose the response in the other final ray so that
the two counts agree after truncation at $T_d$. Each ray allows every
finite offset and remains infinite after the move.

The initial $T_d$-truncations give the invariant before the first
move, so Duplicator wins the $d$-round game.

For (\ref{item:ord-types}), agreement at all truncations gives
agreement at every quantifier rank by (\ref{item:ord-rank}).
Conversely, the requirement that a given gap contain at least $L$
points is first-order expressible.
For (\ref{item:ord-equivalence}), the one-point tuples of least
elements have initial gap zero, so (\ref{item:ord-types}) gives
equality of their complete types. In particular the orders agree
on all sentences.
\end{proof}

Lupini uses the indexed coefficients, the order of their supports,
and the Hilbert Gram matrices to determine tuple types in countably
saturated models of his densely ordered expansion
\cite[Lemma~9.1]{LupiniNormFormulas2026}.

\begin{lemma}[Matching indexed coefficients and Hilbert Gram matrices]
\label{lem:orbit}
Let $\bar a,\bar b$ be finite Hilbert tuples in an elementary
extension $\mathcal M_0$ of $\mathcal M$. Suppose their full Hilbert Gram matrices
agree and there is a partial elementary map of the pure index
order with domain $J(\bar a)$ and range $J(\bar b)$, given by a
bijection $f:J(\bar a)\to J(\bar b)$. Assume that for every
$\alpha\in J(\bar a)$ and every $r$,
\begin{equation}\label{eq:coeffmap}
 \langle a_r,E_\alpha\rangle
   =\langle b_r,E_{f(\alpha)}\rangle,\qquad
 \langle a_r,F_\alpha\rangle
   =\langle b_r,F_{f(\alpha)}\rangle.
\end{equation}
Then $\bar a$ and $\bar b$ have the same complete type in $\mathcal M_0$.
They also have the same complete Banach type in any choice of Banach
ball sorts containing the corresponding components of both tuples.
\end{lemma}
\begin{proof}
Let $I_0$ be the index order of $\mathcal M_0$. For every $\alpha\in J(\bar a)$,
introduce a constant interpreted as $\alpha$ in one copy of $I_0$
and as $f(\alpha)$ in another. The two expanded orders are elementarily
equivalent because every sentence mentions only finitely many constants
and $f$ preserves every pure-order formula on finite tuples from its
domain. By the Keisler--Shelah theorem \cite{Shelah1971}, there is an
ultrafilter $\mathcal V$ such that these expanded orders have
isomorphic ultrapowers by $\mathcal V$. The isomorphism of the
underlying orders is an automorphism $\sigma$ of $I_0^{\mathcal V}$
sending the diagonal image of each $\alpha\in J(\bar a)$ to that
of $f(\alpha)$. Form $\mathcal M_0^{\mathcal V}$ using the same ultrafilter,
so its index order is $I_0^{\mathcal V}$, and identify vectors of
$\mathcal M_0$ with their diagonal images.

To compute the indexed support of a diagonal vector $a$, take
$\beta=[\beta_s]_{\mathcal V}\in I_0^{\mathcal V}$.
If $\langle a,E_\beta\rangle\ne0$, fix
$0<t<|\langle a,E_\beta\rangle|$. By \L o\'s's theorem for continuous
logic \cite[Theorem~5.4]{BBHU08},
$|\langle a,E_{\beta_s}\rangle|>t$ on a set in $\mathcal V$.
Bessel's inequality gives at most $\h(a)^2/t^2$ indices in $I_0$
whose $E$-coefficients have magnitude greater than $t$.
Thus $\beta_s$ takes one fixed value on a set in $\mathcal V$,
so $\beta$ is diagonal. The same reasoning applies if
$\langle a,F_\beta\rangle\ne0$. Coefficients
at diagonal indices retain their original values. Thus the
indexed supports of $\bar a$ and $\bar b$ in the ultrapower
are the diagonal images of $J(\bar a)$ and $J(\bar b)$.

In the Hilbert sort of $\mathcal M_0^{\mathcal V}$, let $A$ be the closed span
of $E_\alpha,F_\alpha$ for $\alpha\in I_0^{\mathcal V}$ and let
$D=A^\perp$, with orthogonal projections $\pi_A$ and $\pi_D$.
The maps $E_\alpha\mapsto E_{\sigma(\alpha)}$ and
$F_\alpha\mapsto F_{\sigma(\alpha)}$ extend to a surjective linear
isometry $W_A:A\to A$. For each $\beta\in I_0^{\mathcal V}$,
\[
 \langle W_A\pi_Aa_r,E_\beta\rangle
 =\langle a_r,E_{\sigma^{-1}(\beta)}\rangle
 =\langle b_r,E_\beta\rangle
 =\langle \pi_Ab_r,E_\beta\rangle.
\]
The middle equality follows from \eqref{eq:coeffmap} on the matched
supports and from vanishing of both coefficients off those supports.
The same equalities hold with $F$ in place of $E$.
Since the indexed vectors span a dense subspace of $A$, it follows that
$W_A\pi_Aa_r=\pi_Ab_r$ for every $r$.

Orthogonality gives
\[
 \langle \pi_Da_r,\pi_Da_s\rangle
 =\langle a_r,a_s\rangle-\langle \pi_Aa_r,\pi_Aa_s\rangle
 =\langle b_r,b_s\rangle-\langle \pi_Ab_r,\pi_Ab_s\rangle
 =\langle \pi_Db_r,\pi_Db_s\rangle.
\]
The equal Gram matrices give equal norms for corresponding linear
combinations, so $\pi_Da_r\mapsto\pi_Db_r$ defines a linear isometry
from the first span onto the second. Let $D_0$ be the span of both
projected tuples. The orthogonal complements of the two spans in
$D_0$ have the same dimension. Mapping an orthonormal basis of the
first complement onto an orthonormal basis of the second extends the
isometry to all of $D_0$. Use the identity on $D\cap D_0^\perp$ to obtain
a surjective linear isometry $W_D:D\to D$.

The direct sum $W=W_A\oplus W_D$ sends $\bar a$ to $\bar b$.
It preserves the Hilbert operations and inner product, and it
intertwines both coordinate maps with $\sigma$. Since $\sigma$
permutes the ordered pairs of indices, the supremum in
\eqref{eq:qstar} gives $q(Wz)=q(z)$ for every $z$. Thus
$(W,\sigma)$ is an automorphism of $\mathcal M_0^{\mathcal V}$. The diagonal
embedding $\mathcal M_0\preccurlyeq \mathcal M_0^{\mathcal V}$
\cite[Corollary~5.5]{BBHU08} gives equality of
complete types in $\mathcal M_0$, and Lemma~\ref{lem:interp} gives equality
of types in the interpreted Banach structure.
\end{proof}

\subsection{Uniform approximation of formula values}
For a finite set $S=\{s_1<\cdots<s_t\}\subseteq\mathbb N$ and
integer $L\ge1$ define the tuple of $L$-truncated gaps, with
the cardinality of $S$ as its first entry,
\begin{equation}\label{eq:truncated-gaps}
 \operatorname{gap}_L(S)=\bigl(
 t,\min\{s_1-1,L\},
 \bigl(\min\{s_{j+1}-s_j-1,L\}\bigr)_{1\le j<t}\bigr)
 \qquad(t\ge1).
\end{equation}
Put $\operatorname{gap}_L(\varnothing)=(0)$. The final ray
does not enter the definition because it is always infinite.

\begin{lemma}[Restriction of truncated gaps]\label{lem:gap-restriction}
Let $S,S'$ be finite subsets of $\mathbb N$ and let $L\ge1$ be an
integer. Suppose $\operatorname{gap}_L(S)=\operatorname{gap}_L(S')$ and
let $f:S\to S'$ be the increasing bijection. For every $T\subseteq S$,
\begin{equation}\label{eq:gap-restriction}
 \operatorname{gap}_L(T)=\operatorname{gap}_L(f(T)).
\end{equation}
\end{lemma}
\begin{proof}
If $T$ is empty, both tuples are $(0)$. Otherwise write the gaps
of $S=\{s_1<\cdots<s_t\}$ as $g_0=s_1-1$ and
$g_j=s_{j+1}-s_j-1$ for $1\le j<t$. If $s_p$ is the first
point of $T$, its initial gap is
$\sum_{j=0}^{p-1}g_j+(p-1)$. If $s_p,s_q$ are successive
points of $T$, the gap between them is
$\sum_{j=p}^{q-1}g_j+(q-p-1)$. In each formula the added
integer counts the points of $S$ that were deleted in that gap.
The corresponding formulas for $f(T)$ have the same summation
ranges and the same numbers of deleted points. If a summand
is at least $L$, the truncated sum is $L$. If all summands
are less than $L$, truncating them leaves the sum unchanged.
Thus the truncated sums depend only on the original truncated
gaps and the numbers of deleted points.
\end{proof}

The indexed coefficients of a sequence can tend uniformly to zero
while its Hilbert norms stay bounded away from zero.
For example, in the Banach ball of radius $2$, take
$a^{(n)}=n^{-1/2}\sum_{i=1}^n E_i$ and $b^{(n)}=0$.
For any fixed positive threshold, all indexed coefficients of
both vectors are smaller in magnitude than the threshold for
sufficiently large $n$.
Nevertheless $\Nrm(a^{(n)})=\sqrt{1+1/n}$ and $\Nrm(b^{(n)})=0$.
In a nonprincipal ultrapower over $\mathbb N$, the two sequences
give vectors with empty indexed supports and
Hilbert norms $1$ and $0$, respectively. The full Hilbert Gram
matrix controls the components orthogonal to the indexed span,
as in Lemma~\ref{lem:orbit}.

\begin{lemma}[Uniform approximation of formula values]\label{lem:formula-approximation}
Let $\varphi(z_1,\ldots,z_m)$ be a Banach-language formula on
fixed ball sorts with radii bounded by $R<\infty$. For every
$\varepsilon>0$, there exist $\delta>0$ and an integer $L\ge1$,
depending only on $(\varphi,R,\varepsilon)$, such that whenever tuples
$\bar a,\bar b$ in the product of the indicated ball sorts in $\Y$
and finite sets $S,S'\subseteq\mathbb N$ satisfy the following conditions
\begin{enumerate}[label=\textup{(\roman*)},leftmargin=3em]
\item every $E$- or $F$-coefficient of $\bar a$ outside $S$,
and of $\bar b$ outside $S'$, has absolute value $<\delta$;
\item $\operatorname{gap}_L(S)=\operatorname{gap}_L(S')$,
with $f:S\to S'$ the increasing bijection;
\item at corresponding indices $\alpha\in S$ and $f(\alpha)\in S'$,
all coefficients of corresponding tuple components along $E$ and
$F$ differ by $<\delta$;
\item all corresponding entries of the \emph{full Hilbert Gram
matrices} of $\bar a$ and $\bar b$ differ by $<\delta$,
\end{enumerate}
we have
\begin{equation}\label{eq:uniformapprox}
 \big|\varphi^{\Y}(\bar a)-\varphi^{\Y}(\bar b)\big|
 <\varepsilon.
\end{equation}
The constants $\delta,L$ are uniform over all finite $S,S'$.
\end{lemma}
\begin{proof}
Suppose the conclusion fails. For each $n\ge1$, take a
counterexample with tuples $\bar a^{(n)},\bar b^{(n)}$ and finite
sets $S_n,S'_n$ satisfying the hypotheses for $\delta=1/n$ and
$L=n$. Let $f_n:S_n\to S'_n$ be the increasing bijection. Thus
\begin{equation}\label{eq:contra}
 |\varphi^{\Y}(\bar a^{(n)})-
   \varphi^{\Y}(\bar b^{(n)})|\ge\varepsilon.
\end{equation}
Let $\cU$ be a nonprincipal ultrafilter on $\mathbb N$ and form
$\mathcal M^{\cU}$. Set $a_r=[a_r^{(n)}]_{\cU}$ and
$b_r=[b_r^{(n)}]_{\cU}$ for $1\le r\le m$. By (iv) and
\L o\'s's theorem, the full Hilbert Gram matrices of $\bar a$
and $\bar b$ agree.

Let $\alpha=[i_n]_{\cU}\in J(\bar a)$. Some coefficient of
$\bar a$ at $\alpha$ is nonzero. By \L o\'s's theorem, the
corresponding coefficients in the factors have absolute value
at least $t$ on a set in $\cU$, for some $t>0$. On this set,
condition (i) forces $i_n\in S_n$ whenever $1/n<t$. Thus
$\{n:i_n\in S_n\}\in\cU$.
Following Lupini \cite[Lemma~9.2]{LupiniNormFormulas2026}, set
\[
 f(\alpha)=[j_n]_{\cU},\qquad
 j_n=\begin{cases}
       f_n(i_n),&i_n\in S_n,\\
       1,&i_n\notin S_n.
     \end{cases}
\]
Two representatives of $\alpha$ agree on a set in $\cU$, and
so do the corresponding sequences $(j_n)$. By (iii) and
\L o\'s's theorem,
\begin{equation}\label{eq:mapultra}
 \langle a_r,E_\alpha\rangle=\langle b_r,E_{f(\alpha)}\rangle,\qquad
 \langle a_r,F_\alpha\rangle=\langle b_r,F_{f(\alpha)}\rangle
 \qquad(r\le m).
\end{equation}
In particular, $f(\alpha)\in J(\bar b)$. The same definition
with $f_n^{-1}$ gives a map from $J(\bar b)$ to $J(\bar a)$.
The two maps are inverse since the factor maps are inverse on
the sets where the representatives belong to $S_n$ or $S'_n$,
respectively.

Fix a pure order formula $\theta(x_1,\ldots,x_k)$ of quantifier
rank $d$ and a tuple $\alpha_1,\ldots,\alpha_k$ from $J(\bar a)$,
with $\alpha_\ell=[i_{\ell,n}]_{\cU}$. There is a set $Q\in\cU$
on which $i_{\ell,n}\in S_n$ for every $\ell$ and $n\ge 2^d-1$.

For $n\in Q$, let $T_n=\{i_{1,n},\ldots,i_{k,n}\}$. By (ii) and
Lemma~\ref{lem:gap-restriction}, the gaps of $T_n$ and $f_n(T_n)$ agree
after truncation at $n$, hence at $2^d-1$. By
Lemma~\ref{lem:ord}(\ref{item:ord-rank}),
their increasing enumerations agree on all order formulas of rank at
most $d$. The tuple $(i_{1,n},\ldots,i_{k,n})$ and its image under
$f_n$ use the same positions in the two enumerations, including any
repetitions. Hence
\[
 \mathbb N\models\theta(i_{1,n},\ldots,i_{k,n})
 \quad\Longleftrightarrow\quad
 \mathbb N\models\theta(f_n(i_{1,n}),\ldots,f_n(i_{k,n})).
\]
By \L o\'s's theorem, $(\alpha_1,\ldots,\alpha_k)$ and
$(f(\alpha_1),\ldots,f(\alpha_k))$ agree on $\theta$.
Thus $f$ is partial elementary in the pure order.

By Lemma~\ref{lem:orbit}, the tuples $\bar a$ and $\bar b$ have
the same complete Banach type. By Lemma~\ref{lem:interp}, the
interpreted Banach structure in $\mathcal M^{\cU}$ is the Banach
ultrapower of $Y$, where \L o\'s's theorem and \eqref{eq:contra}
give $|\varphi(\bar a)-\varphi(\bar b)|\ge\varepsilon$, a contradiction.
\end{proof}

\subsection{Polynomial bounds for finite sets and Hilbert inner products}
\subsubsection{Relative positions and truncated gaps}
\begin{lemma}\label{lem:relative-positions}
Fix integers $k_0,L\ge1$ and a finite $B\subseteq\mathbb N$.
Partition the sets $T\subseteq\mathbb N$ with
$|T|\le k_0$ by putting two sets in the same class if they agree on
the following three quantities.
\begin{enumerate}[label=\textup{(\alph*)},leftmargin=3em]
\item The intersection $T\cap B$.
\item The number of points of $T$ in each complementary interval of $B$.
\item The tuple $\operatorname{gap}_L(B\cup T)$.
\end{enumerate}
There are at most $C_{k_0,L}(|B|+1)^{k_0}$ classes, with $C_{k_0,L}$
independent of $B$. If $T,T'$ belong to the same class, the
increasing bijection
$\pi:B\cup T\to B\cup T'$ fixes each point of $B$, maps $T$
onto $T'$, and for \emph{every} $S\subseteq B$ its restriction
$S\cup T\to S\cup T'$ preserves the $L$-truncated initial and
successive gaps.
\end{lemma}
\begin{proof}
In increasing order, each point of $T$ belongs to one of the $|B|$
singletons of $B$ or to one of its $|B|+1$ complementary intervals.
For $|T|=t$, there are at most $(2|B|+1)^t$ such lists. Each list
determines (a) and (b). For fixed (a) and (b), the cardinality
of $B\cup T$ is fixed. If a complementary interval contains
$u\ge1$ points of $T$, there are $u-1$ gaps between those points
and at most one further gap at each end of the interval.
In the initial interval we count the gap before the first point,
and in the final ray we omit the gap after the last point.
Over all complementary intervals containing points of $T$,
at most $2k_0$ gap lengths can vary, since $u+1\le2u$.
Each truncated gap length has at most $L+1$ possible values.
The gaps in intervals containing no point of $T$ are determined
by $B$.
Thus the total is bounded by
\[
 (k_0+1)(2|B|+1)^{k_0}(L+1)^{2k_0}
 \le C_{k_0,L}(|B|+1)^{k_0}.
\]
By (b), each complementary interval of $B$ contains the
same number of points of $T$ and of $T'$. Thus the increasing
bijection $\pi:B\cup T\to B\cup T'$ fixes $B$ pointwise and
matches the new points in their interval order. By (a), the
members of $B$ that lie in $T$ are exactly those that lie in
$T'$, so $\pi(T)=T'$. Given $S\subseteq B$, we have
$\pi(S\cup T)=S\cup T'$. By (c), the full unions have equal
$L$-truncated gaps. Applying Lemma~\ref{lem:gap-restriction} to
$\pi$ and the subset $S\cup T$ gives equal $L$-truncated gaps
for $S\cup T$ and $S\cup T'$.
\end{proof}

\subsubsection{Simultaneous approximation of Hilbert inner products}
The following Hilbert space estimate is a special case of the
polynomial covering bounds of Zhang \cite[Theorem~4 with $p=2$]{Zhang2002}.
See also Lupini \cite[Lemma~9.3]{LupiniNormFormulas2026}.

\begin{lemma}[Polynomial covering bound for Hilbert inner products]\label{lem:hilbert-net}
Let $R,\rho>0$ and $u_1,\dots,u_{n_0}$ be vectors in an arbitrary
real Hilbert space $H_0$, with $\|u_i\|\le R$.
Then
\[
 \mathcal C=\{(\langle u_1,v\rangle,\dots,
 \langle u_{n_0},v\rangle):\|v\|\le R\}\subseteq\mathbb R^{n_0}
\]
has a radius-$\rho$ net in $\ell_\infty^{n_0}$ with at most
$(n_0+2)^{C_{R,\rho}}$ centers. The constant $C_{R,\rho}$
depends only on $R,\rho$, not on $n_0$ or $\dim H_0$.
The centers can be chosen from $\mathcal C$ itself.
\end{lemma}
\begin{proof}
For a vector $v$ of norm at most $R$, we seek a span $W$ of a few
of the $u_i$ such that the residual $v-\pi_Wv$ has inner product
at most $\rho/2$ in absolute value with every $u_i$. Here $\pi_W$
denotes the orthogonal projection. Start with $W=\{0\}$.
If $r=v-\pi_Wv$ satisfies $|\langle u_i,r\rangle|>\rho/2$ for some
$i$, enlarge $W$ by $u_i$. The new direction orthogonal to $W$
is $w=u_i-\pi_Wu_i$. Since $r\perp W$, we have
$\langle r,w\rangle=\langle r,u_i\rangle$ and $0<\|w\|\le R$.
Projecting onto the enlarged span decreases the squared residual
norm by
\[
 \frac{|\langle r,w\rangle|^2}{\|w\|^2}>\frac{\rho^2}{4R^2}.
\]
The initial squared norm is at most $R^2$, so the procedure
stops after at most $m=\lceil4R^4/\rho^2\rceil$ enlargements.

There are at most $(n_0+1)^m$ spans of subsets of at most $m$
of the $u_i$, by listing the indices in increasing order and
padding each list with zeros. In the radius-$R$ ball of each
such span $W$, take a maximal set of points at pairwise distance
greater than $\rho/(2R)$. The balls of radius $\rho/(4R)$ around
these points are disjoint and lie in the radius-$(R+\rho/(4R))$
ball of $W$. Comparing volumes bounds the number of points by
$(1+4R^2/\rho)^m$. For $W=\{0\}$ use the singleton $\{0\}$.
Maximality gives a $\rho/(2R)$-net of the ball of $W$.

Use the evaluations of these net points as centers in $\mathcal C$.
For any $v$, take its final span $W$ and a net point $z$ within
$\rho/(2R)$ of $\pi_Wv$, which has norm at most $R$. For every $i$,
\[
 |\langle u_i,v-z\rangle|
 \le |\langle u_i,v-\pi_Wv\rangle|+\|u_i\|\,\|\pi_Wv-z\|
 \le \rho .
\]
The number of centers is at most $(n_0+1)^m(1+4R^2/\rho)^m$,
which is bounded by $(n_0+2)^{C_{R,\rho}}$ for a suitable
$C_{R,\rho}$.
\end{proof}

\subsection{Subsets determined between two thresholds}

Fix a Banach-language formula $\varphi(\bar x;\bar y)$ on fixed ball
sorts. For tuples $\bar a_1,\ldots,\bar a_n$ in $\Y$, their evaluation
set consists of the vectors
$(\varphi(\bar a_i;\bar b))_{i=1}^n$ as $\bar b$ ranges over its
ball sorts in $\Y$. For $r<s$, this set determines
$U\subseteq\{1,\ldots,n\}$ between $r$ and $s$ in the sense of
Ben Yaacov \cite[Definition~2.3]{BenYaacovVC} if some $\bar b$ satisfies
\begin{equation}\label{eq:IP}
 \varphi(\bar a_i;\bar b)\ge s\quad(i\in U),\qquad
 \varphi(\bar a_i;\bar b)\le r\quad(i\notin U).
\end{equation}

\begin{proposition}\label{prop:determined-subsets}
Fix $\varphi(\bar x;\bar y)$, its ball sorts, and $r<s$. There are
constants $C,A$ depending only on these choices such that, for any
tuples $\bar a_1,\ldots,\bar a_n$ in $\Y$, their evaluation set determines
at most $C(n+1)^A$ subsets between $r$ and $s$.
\end{proposition}

\begin{proof}
Let $m$ and $k$ be the lengths of $\bar x$ and $\bar y$, and let
$R\ge1$ bound the radii of their Banach ball sorts. If $m=0$,
each $\bar a_i$ is the empty tuple, so only $\varnothing$ and
$\{1,\ldots,n\}$ can be determined. If $k=0$, the evaluation set
consists of one vector and determines at most one subset. Assume $m,k\ge1$.

Put $g=s-r$. Apply Lemma~\ref{lem:formula-approximation} to
$(\bar x,\bar y)$ with error $g/4$ and ball bound $R$ to obtain
$\delta>0$ and $L\ge1$. Set $\eta=\delta/4$. For a tuple $\bar z$,
let $S(\bar z)$ consist of the indices at which some $E$- or
$F$-coefficient of a component has absolute value at least $\eta$.
Put
\[
 S_i=S(\bar a_i),\qquad S_0=\bigcup_{i=1}^n S_i,\qquad
 T_{\bar b}=S(\bar b).
\]
By Bessel's inequality and $\h\le\Nrm$, we have
$|S_i|\le mR^2/\eta^2$ and $|T_{\bar b}|\le k_0$, where
$k_0=\lceil kR^2/\eta^2\rceil$.

We partition the parameter tuples so that the values of
$\varphi(\bar a_i;\cdot)$ differ by less than $g/4$ within each class
for every $i$. For each $\bar b$ we use three quantities, as in
Lupini \cite[Theorem~9.4]{LupiniNormFormulas2026}.
\begin{enumerate}[label=(\alph*),leftmargin=3em]
\item The set $T_{\bar b}$ belongs to a class in the partition from
Lemma~\ref{lem:relative-positions} with $B=S_0$, size bound $k_0$,
and gap cutoff $L$.
\item Fix finite partitions of $[-R,R]$ and $[-R^2,R^2]$ into intervals
of length at most $\eta$. The coefficients on $T_{\bar b}$ give a
list of intervals in the first partition. The indices occur in
increasing order and the components in tuple order, with the
$E$-coefficient before the $F$-coefficient for each component.
The entries of the Hilbert Gram matrix give a list of intervals in the
second partition, ordered by the pairs of component indices.
\item For the $nm$ vectors $a_{i,\ell}$, fix one enumerated
radius-$\eta$ net from Lemma~\ref{lem:hilbert-net} for all parameter
tuples. Each $b_t$ gives the first center within $\eta$ of
$(\langle a_{i,\ell},b_t\rangle)_{i,\ell}$.
\end{enumerate}
Two parameter tuples belong to the same class exactly when the class
in (a), the lists in (b), and the ordered $k$-tuple of centers in (c)
agree. The number of possibilities in (b) is bounded independently of $n$ and
the tuples $\bar a_i$, since there are at most $2kk_0$ coefficients and $k^2$
Gram entries. By Lemmas~\ref{lem:relative-positions}
and~\ref{lem:hilbert-net}, the number of classes in (a) and the number
of $k$-tuples of net centers in (c)
are at most $C_{k_0,L}(|S_0|+1)^{k_0}$ and $(nm+2)^{kC_{R,\eta}}$,
respectively. Since $|S_0|\le nmR^2/\eta^2$, there are at most
$C(n+1)^A$ classes, with $C,A$ independent of $n$ and the tuples $\bar a_i$.

Take $\bar b,\bar b'$ in the same class. By (a) and
Lemma~\ref{lem:relative-positions}, the increasing bijection
\[
 \pi:S_0\cup T_{\bar b}\longrightarrow S_0\cup T_{\bar b'}
\]
fixes $S_0$ pointwise and maps $T_{\bar b}$ onto $T_{\bar b'}$.
The two unions have equal $L$-truncated gaps. For every $i$ we use
these same two sets and the same map $\pi$ to compare
$(\bar a_i,\bar b)$ with $(\bar a_i,\bar b')$ in
Lemma~\ref{lem:formula-approximation}.

Outside the respective unions all coefficients of the combined tuples
have absolute value less than $\eta$, since $S_i\subseteq S_0$.
At matched indices in $S_0$ the coefficients of $\bar a_i$ agree exactly.
At matched indices outside $S_0$ both coefficients of $\bar a_i$ have
absolute value less than $\eta$. By (b), the parameter coefficients
at matched indices in $T_{\bar b}$ and $T_{\bar b'}$ differ by at
most $\eta$. Outside these respective sets both have absolute value
less than $\eta$. Thus every matched coefficient difference is at
most $2\eta$. In the full Hilbert
Gram matrices, the entries among components of $\bar a_i$ agree, and
the entries among the parameter vectors differ by at most $\eta$
by (b). The mixed entries differ by at most $2\eta$ by (c), since
the two vectors of mixed inner products for each component are
within $\eta$ of the same center.

Since $2\eta<\delta$, Lemma~\ref{lem:formula-approximation} gives
\[
 |\varphi(\bar a_i;\bar b)-\varphi(\bar a_i;\bar b')|<g/4
 \qquad(1\le i\le n).
\]
If $\bar b$ and $\bar b'$ satisfy \eqref{eq:IP} for distinct subsets,
their formula values differ by at least $g$ at any index belonging to
exactly one of the two subsets. Thus each parameter class determines
at most one subset between $r$ and $s$.
\end{proof}

\subsection{NIP of the Banach theory}

\begin{theorem}\label{thm:main}
The complete Banach theory $\operatorname{Th}(\Y)$ is NIP.
\end{theorem}

\begin{proof}
Suppose a Banach-language formula $\varphi(\bar x;\bar y)$ has IP in
$\operatorname{Th}(\Y)$ with thresholds $r<s$. Fix $0<\gamma<(s-r)/2$.
By \cite[Lemma~2.3(ii)]{PoitevinNormFormula2026v2}, for each $n$
there are tuples $\bar a_1,\ldots,\bar a_n$ in $\Y$ whose evaluation
set determines every subset of $\{1,\ldots,n\}$ between $r+\gamma$
and $s-\gamma$. Proposition~\ref{prop:determined-subsets} then gives
$2^n\le C(n+1)^A$, with $C,A$ independent of $n$. This fails for
sufficiently large $n$.
\end{proof}

\section{Separable models of \texorpdfstring{$\operatorname{Th}(Y)$}{Th(Y)}}
\subsection{Recovering Hilbert geometry and the index order}
\label{sec:interpret-order}
We now recover the Hilbert norm, the graph on ordered pairs of indices,
and the index order from the Banach norm and the vector space
operations.
Throughout this section, $B_R$ is the closed ball of radius $R$ for
$\Nrm$; all definability assertions refer to the pure Banach language.

Let $I$ be a countable model of $\operatorname{Th}(\mathbb N,<)$, with
least point $c_I$. Write $E_i,F_j$ for the coordinate vectors
in the first and second copies of $\ell_2(I)$ and $\h_I$ for the
Hilbert norm on $H_I=\ell_2(I)\oplus_2\ell_2(I)$. The seminorm and
Banach norm on $H_I$ are
\[
 q_I(u,v)=\sup_{i<_I j}|u_i+v_j|,\qquad
 \Nrm_I(x)=\bigl(\h_I(x)^2+q_I(x)^2\bigr)^{1/2}.
\]
Write $Y_I=(H_I,\Nrm_I)$, so $Y_{\mathbb N}=\Y$, and let
\[
 Y_I^\circ=\{(u,v)\in Y_I:v_{c_I}=0\}.
\]
No pair $i<_I j$ has $j=c_I$, so $q_I$ is unchanged by adding a
multiple of $F_{c_I}$. Thus $Y_I\cong Y_I^\circ\oplus_2\mathbb R$.

\begin{lemma}[Intrinsic recovery of $\h$ and $q$]\label{lem:intrinsic-h}
Let $K$ be a real Hilbert space. In $Y_I^\circ\oplus_2K$, write
$h_0$ for the Hilbert sum norm, so
$\Nrm(x)^2=h_0(x)^2+q_I(u,v)^2$ for $x=(u,v,k)$.
For every positive rational $R$, let $B_R$ be the closed ball
of radius $R$ for the Banach norm on this sum, and put $t_R=1/(4R)$.
On $B_R$ we have
\begin{equation}\label{eq:intrinsic-h}
 h_0(x)^2=\inf_{z\in B_3}
 \frac{\Nrm(z+t_Rx)^2+\Nrm(z-t_Rx)^2-2\Nrm(z)^2}{2t_R^2}
\end{equation}
and $q_I(u,v)=\sqrt{\Nrm(x)^2-h_0(x)^2}$.
In particular the Hilbert norm $\h$, the Hilbert inner product,
and the seminorm $q$ of $\Y$ are uniformly definable on bounded sorts
of $\operatorname{Th}(\Y)$. The Hilbert identities and the norm
comparisons also hold in every model.
\end{lemma}
\begin{proof}
By the parallelogram identity, the contribution of $h_0^2$ to the
numerator in \eqref{eq:intrinsic-h} is $2t_R^2h_0(x)^2$. The
contribution of $q_I^2$ is nonnegative by convexity, so the quotient
is bounded below by $h_0(x)^2$.

For the reverse inequality, let $i_1=c_I$ and let $i_{n+1}$ be the
immediate successor of $i_n$. Set $z_n=(E_{i_n}+F_{i_{n+1}},0)$.
Then $h_0(z_n)=\sqrt2$ and the coordinate sum for $(i_n,i_{n+1})$
is $2$. All other sums have absolute value at most $1$, so
$\Nrm(z_n)=\sqrt6<3$. For $x=(u,v,k)\in B_R$, adding or subtracting
$t_R x$ changes each coordinate sum by at most
$t_R q_I(u,v)\le t_R\Nrm(x)\le1/4$. Thus the sum for $(i_n,i_{n+1})$
is at least $7/4$, while all other sums have absolute value at most
$5/4$. The seminorm of the coordinate part of $z_n\pm t_R x$ is
therefore $2\pm t_R(u_{i_n}+v_{i_{n+1}})$. The quotient at $z_n$ is
\[
 h_0(x)^2+(u_{i_n}+v_{i_{n+1}})^2.
\]
Both coefficients tend to zero by square summability, so the
infimum equals $h_0(x)^2$.

For $I=\mathbb N$ and $K=\mathbb R$, the identity applies to $\Y$.
For every $z\in B_3$ and $x\in B_R$, the terms $z\pm t_Rx$ belong to
$B_{13/4}$, so all quantification and connectives occur on fixed
bounded sorts. The right-hand side of \eqref{eq:intrinsic-h} defines
$\h^2$ on $\Y$. Taking a square root gives $\h$, and polarization
gives $\langle\cdot,\cdot\rangle_{\h}$. Finally
$q=\sqrt{(\Nrm^2-\h^2)_+}$ is uniformly definable on bounded
sorts. The norm axioms, the parallelogram identity for $\h$,
the seminorm axioms for $q$, and the norm comparisons transfer
as identities and inequalities on bounded sorts. In every model,
the Jordan--von Neumann theorem \cite[Theorem~I]{JordanVonNeumann1935}
makes $\h$ an inner product norm. It is complete because it is
equivalent to the Banach norm.
\end{proof}

\begin{lemma}[Definability of the signed pair vectors]\label{lem:edge-set}
For $i<j$ put
\[
 e_{ij}=\frac{E_i+F_j}{\sqrt2},\qquad
 \mathcal E=\{\pm e_{ij}:i<j\}\subseteq B_2.
\]
The function
\begin{equation}\label{eq:edge-distance}
 d_{\mathcal E}(x)=\sqrt{\bigl(\h(x)^2+1-\sqrt2\,q(x)\bigr)_+}
 \quad(x\in B_2)
\end{equation}
is the Hilbert distance from $x$ to $\mathcal E$.  In every model
$M\models\operatorname{Th}(\Y)$ the zero set
$\mathcal E(M)=d_{\mathcal E}^{-1}(0)$ is a definable set. Its
distance predicate and restricted quantifiers are definable
uniformly in $M$.
\end{lemma}
\begin{proof}
For $x=(u,v)\in\Y$ and $e\in\mathcal E$, we have
$\h(x-e)^2=\h(x)^2+1-2\langle x,e\rangle_{\h}$.
For each $i<j$ the two signs of $e_{ij}$ give inner products
$\pm(u_i+v_j)/\sqrt2$. Thus
$\sup_{e\in\mathcal E}\langle x,e\rangle_{\h}=q(x)/\sqrt2$ and
\[
 \inf_{e\in\mathcal E}\h(x-e)^2
 =\h(x)^2+1-\sqrt2\,q(x).
\]
Distinct signed pair vectors have Hilbert distance at least one.
Thus $\mathcal E$ is closed, and the distance formula vanishes
exactly on $\mathcal E$.
In every model, $q\le\sqrt2\h$, so the radicand is at least
$(\h(x)-1)^2$. In particular every member of $\mathcal E(M)$
has Hilbert norm one.

In $\Y$, the metric-distance identities
\begin{align}
 |d_{\mathcal E}(x)-d_{\mathcal E}(y)|&\le\h(x-y),\label{eq:edge-Lipschitz}\\
 \inf_{z\in B_2}\max\bigl\{d_{\mathcal E}(z),\,
       |\h(x-z)-d_{\mathcal E}(x)|\bigr\}&=0 \qquad(x\in B_2)
 \label{eq:edge-approx-near}
\end{align}
hold. Both transfer to every model because $d_{\mathcal E}$ and $\h$ are
uniformly definable. Fix $x\in B_2$ in such a model $M$ and
$\varepsilon>0$. Choose $\eta_n>0$ for $n\ge0$ with
$2\sum_{n\ge0}\eta_n<\varepsilon$. Starting with $x_0=x$, apply
\eqref{eq:edge-approx-near} to each $x_n$ to obtain
$x_{n+1}\in B_2$ with $d_{\mathcal E}(x_{n+1})<\eta_n$ and
$\h(x_n-x_{n+1})<d_{\mathcal E}(x_n)+\eta_n$. In particular, for $n\ge1$
the last bound is less than $\eta_{n-1}+\eta_n$. Thus
\[
 \sum_{n=0}^{\infty}\h(x_n-x_{n+1})
 \le d_{\mathcal E}(x)+2\sum_{n=0}^{\infty}\eta_n
 <d_{\mathcal E}(x)+\varepsilon.
\]
For $m>n$, the triangle inequality bounds $\h(x_m-x_n)$
by the tail of this series starting at $n$. The sequence is
therefore Cauchy for $\h$ and for the equivalent norm $\Nrm$.
Its limit $e$ lies in $B_2$. Since $d_{\mathcal E}(x_{n+1})<\eta_n\to0$,
\eqref{eq:edge-Lipschitz} gives $d_{\mathcal E}(e)=0$. The bound on the series also gives
$\h(x-e)<d_{\mathcal E}(x)+\varepsilon$. Since $\varepsilon$ was arbitrary,
$\dist_\h(x,\mathcal E(M))\le d_{\mathcal E}(x)$. For every $e'\in\mathcal E(M)$,
\eqref{eq:edge-Lipschitz} gives $d_{\mathcal E}(x)\le\h(x-e')$, proving equality.

If $\psi(z,\bar y)$ is a $[0,1]$-valued definable predicate on
fixed sorts, define
\[
 a_n(\bar y)=\inf_{z\in B_2}
 \min\bigl\{1,\psi(z,\bar y)+n\,d_{\mathcal E}(z)\bigr\}.
\]
Let $\omega$ be a modulus of continuity for $\psi$ in its first
variable for the Banach metric. Evaluating the expression at signed
vectors gives $a_n(\bar y)\le\inf_{e\in\mathcal E(M)}\psi(e,\bar y)$.
If $d_{\mathcal E}(z)\ge1/n$, the term is one. If
$d_{\mathcal E}(z)<1/n$, there is an $e\in\mathcal E(M)$ with
$\h(z-e)<1/n$, so $\Nrm(z-e)<\sqrt3/n$. In that case the term
is at least $\psi(z,\bar y)\ge\psi(e,\bar y)-\omega(\sqrt3/n)$.
Taking the infimum over $z$ gives the uniform estimate
\[
 0\le \inf_{e\in\mathcal E(M)}\psi(e,\bar y)-a_n(\bar y)
 \le\omega(\sqrt3/n).
\]
Complementation gives restricted suprema. Taking $\psi$ to
be a rescaled Banach distance gives the Banach distance predicate.
\end{proof}

\begin{lemma}[First-row order and recovery of the two indices]
\label{lem:graph-indices}
Let $I$ be a linear order with least point $c$ and an immediate
successor $s(i)$ for every $i$. Regard each pair $(i,j)$ with $i<j$
as a cell in row $i$ and column $j$. Write $e\parallel f$ when $e$
and $f$ share a row or a column. This relation is reflexive.
In the language with the single relation $\parallel$, the first row
$P=\{p_j=(c,j):j>c\}$ and its order given by
$p_j\prec p_k$ if and only if $j<k$ are definable without parameters.

For $e=(i,j)$ outside $P$, the unique $p\in P$ satisfying $e\parallel p$
is $p_j$. The $\prec$-least $p\in P$ for which some $g\notin P$
satisfies $g\parallel e$ and $g\parallel p$ is $p_{s(i)}$. These two first-row
cells are definable from $e$ and determine $e$ uniquely among the
cells outside $P$.
\end{lemma}
\begin{proof}
The following formulas quantify over cells.
\begin{align*}
 \theta_0(o)&\ \Longleftrightarrow\
 \forall u\,\forall v\,\bigl((u\parallel o\wedge v\parallel o)
                                      \Rightarrow u\parallel v\bigr),\\
 \theta_P(p)&\ \Longleftrightarrow\
 \exists o\,\bigl(\theta_0(o)\wedge p\parallel o\bigr).
\end{align*}
Column $s(c)$ contains only $(c,s(c))$, so every cell related to
$(c,s(c))$ lies in row $c$. Hence $\theta_0$ holds at $(c,s(c))$.
Every other cell $(i,j)$ has $j>s(c)$. Let $k=s(c)$ if $i=c$ and
$k=c$ if $i>c$. Both $(i,s(j))$ and $(k,j)$ share a row or a
column with $(i,j)$. They have different rows and different columns
from each other, so they witness the failure of $\theta_0$ at $(i,j)$.
Thus $\theta_0$ defines the single cell $(c,s(c))$, and $\theta_P$
defines the first row. Write $p\in P$ for $\theta_P(p)$.

For $p\in P$, let $C(p)=\{g\notin P:g\parallel p\}$, the cells in the
column of $p$ below the first row. Define
\[
 \begin{aligned}
 \theta_\prec(p,p')\ \Longleftrightarrow\quad
   &p\in P\wedge p'\in P\wedge p\ne p'\\
   &{}\wedge\forall g\in C(p)\ \exists g'\in C(p')\ (g\parallel g').
 \end{aligned}
\]
For distinct $p_j,p_k\in P$, two cells in their respective columns
are related by $\parallel$ exactly when their row indices coincide. Thus the
quantified condition for $\theta_\prec(p_j,p_k)$ is the inclusion
$\{i:c<i<j\}\subseteq\{i:c<i<k\}$. If $j<k$, the inclusion holds.
If $k<j$, row $k$ belongs to the set for column $j$ and is absent from
the set for column $k$. Hence $\theta_\prec$ defines $\prec$.

For $e=(i,j)$ outside $P$, the cell $p_k$ shares a row or a column
with $e$ exactly when $k=j$. The columns containing a cell in row
$i$ are given by the first-row cells in the set
\[
 A(e)=\{p\in P:\exists g\notin P\,(g\parallel e\wedge g\parallel p)\}.
\]
Indeed, a witness $g$ for $p_k\in A(e)$ must have column $k$.
If $k\ne j$, it must have row $i$ to satisfy $g\parallel e$, so such a
witness exists exactly when $i<k$. If $k=j$, then $i<k$ and $e$
itself is a witness by reflexivity. Thus $A(e)$ consists
of the $p_k$ with $k>i$, and its $\prec$-least element is $p_{s(i)}$.
The two columns $j$ and $s(i)$ determine $e$, since $i$ is the
immediate predecessor of $s(i)$.
\end{proof}

The table summarizes the formulas in Lemma~\ref{lem:graph-indices}.
In the last three rows, $e=(i,j)$ lies outside $P$.
\begin{center}
\begin{tabular}{@{}ll@{}}
Formula & Condition on the indices\\
\hline
$\theta_0(o)$ & $o=(c,s(c))$\\
$\theta_P(p)$ & $p=p_j$ for some $j>c$\\
$g\in C(p_j)$ & $g=(i,j)$ for some $c<i<j$\\
$\theta_\prec(p_j,p_k)$ & $j<k$\\
$p\in P$ and $e\parallel p$ & $p=p_j$\\
$p_k\in A(e)$ & $i<k$\\
$p=\min_\prec A(e)$ & $p=p_{s(i)}$
\end{tabular}
\end{center}

\begin{theorem}[Parameter-free interpretation of the index order]
\label{thm:order-interpretation}
The pure Banach structure $\Y$ interprets $(\mathbb N,<)$ without
parameters on a definable subset of the discrete quotient
$\mathcal E/\{\pm1\}$. In every model of $\operatorname{Th}(\Y)$ the same
interpretation yields a model of $\operatorname{Th}(\mathbb N,<)$.
Every elementary Banach embedding induces elementary embeddings of
the quotient structures and of the interpreted orders.
\end{theorem}
\begin{proof}
For $e,f\in\mathcal E$ define
\[
 \Delta(e,f)=\min\{1,\h(e-f),\h(e+f)\}.
\]
This is the truncated distance between the sign classes, so it is a
pseudometric whose zero classes are the pairs $\{e,-e\}$. Distinct
classes have distance one. The sign-invariant predicate
\[
 \mu([e],[f])=\min\{1,2|\langle e,f\rangle_{\h}|\}
\]
defines the relation $[e]\parallel[f]$ on the quotient
$G=\mathcal E/\{\pm1\}$ by the condition $\mu([e],[f])=1$. In $\Y$ the class
$[e_{ij}]$ corresponds to the pair $(i,j)$. The inner products
\[
 \langle e_{ij},e_{k\ell}\rangle_{\h}
    =\frac12\bigl(\ind_{i=k}+\ind_{j=\ell}\bigr)
\]
show that $\mu$ is $1$ exactly when the two pairs have the same row
or the same column, and is $0$ otherwise. The conditions that
$\Delta$ and $\mu$ take values in $\{0,1\}$ transfer to every model
by the restricted quantifiers of Lemma~\ref{lem:edge-set}.
By Lemma~\ref{lem:graph-indices}, the order on the first row $P$ of
$G$ in $\Y$ is isomorphic to $(\mathbb N,<)$ by $n\mapsto p_{n+1}$.

Classical formulas in the language of $\parallel$ translate to Banach
predicates with $1$ for truth and $0$ for falsehood. The atomic
predicates are $\mu([e],[f])$ for $[e]\parallel[f]$ and $1-\Delta(e,f)$
for equality of the two classes. Negation is $1-a$, conjunction
is $\min$, and implication is $\max\{1-a,b\}$. Existential and
universal quantifiers become respectively $\sup$ and $\inf$ over $\mathcal E$, using
the uniformly definable restricted quantifiers of
Lemma~\ref{lem:edge-set}. An order formula translates to a formula
on $G$ by replacing each instance of the order relation with
$\theta_\prec$ and restricting each quantifier to the first row.
For an existential quantifier the restriction is a conjunction
with $\theta_P$. For a universal quantifier $\theta_P$ is the
antecedent of an implication. Neither defining formula uses parameters.
Thus the interpreted order satisfies $\operatorname{Th}(\mathbb N,<)$
in every model. An elementary Banach embedding $M\to N$ preserves
the values of the Banach predicates translating formulas on $G$,
so it induces an elementary embedding $G(M)\to G(N)$.
Restricting this embedding to the definable
first row gives an elementary embedding of the interpreted orders.
\end{proof}

\begin{corollary}[Uniform definability of the first row]
\label{cor:first-row-definable}
The set $\widetilde P=\{\pm e_{1j}:j\ge2\}\subseteq B_2$ of signed
first-row vectors is uniformly definable in the pure Banach language.
In every model of $\operatorname{Th}(\Y)$ its Hilbert-distance predicate is
uniformly definable on each Banach ball sort $B_R$, and
quantifiers restricted to the set are uniformly definable.
\end{corollary}
\begin{proof}
In every model, let $\widetilde P$ be the preimage of the
first-row subset $P$ of the discrete quotient $G$. The formula
defining $P$ gives a $\{0,1\}$-valued predicate
$\chi_{\widetilde P}$ on $\mathcal E$, with value $1$ exactly
on $\widetilde P$. Fix a Banach ball sort $B_R$. The first row
is nonempty in every model, and for $x\in B_R$ every signed
vector $e$ satisfies $\h(x-e)\le R+1$, since $\h(x)\le\Nrm(x)$
and $\h(e)=1$. Therefore the Hilbert distance from $x$ to
$\widetilde P$ is
\[
 \inf_{e\in\mathcal E}
 \bigl(\h(x-e)+(R+2)(1-\chi_{\widetilde P}(e))\bigr).
\]
For $e\in\widetilde P$ the expression is $\h(x-e)\le R+1$,
and outside $\widetilde P$ it is at least $R+2$.
The infimum is uniformly definable on $B_R$ by
Lemma~\ref{lem:edge-set}, since $e$ ranges over $\mathcal E\subseteq B_2$
and $x-e$ lies in $B_{R+2}$.

For a $[0,1]$-valued definable predicate $\psi(z,\bar y)$ with
$z\in B_2$, we have
\[
 \inf_{z\in\widetilde P}\psi(z,\bar y)
 =\inf_{e\in\mathcal E}
   \min\{1,\psi(e,\bar y)+1-\chi_{\widetilde P}(e)\}.
\]
The term is $\psi(e,\bar y)$ on $\widetilde P$ and $1$ elsewhere
on $\mathcal E$. This gives the restricted infimum by
Lemma~\ref{lem:edge-set}. Complementation gives restricted
suprema, and taking $\psi$ to be a rescaled Banach distance
gives the Banach distance predicate.
\end{proof}

\subsection{Classification and separable spectrum}
\label{sec:separable-classification}
\begin{lemma}[Relabeling coordinates in ultrapowers]\label{lem:ultra-relabel}
Let $I,J$ be countable models of $\operatorname{Th}(\mathbb N,<)$.
For an ultrafilter $\cU$, let $\sigma:I^{\cU}\to J^{\cU}$ be an
order isomorphism.
Relabeling the indexed vectors by $\sigma$ extends to a surjective
linear isometry $Y_I^{\cU}\to Y_J^{\cU}$. The same holds for
$(Y_I^\circ)^{\cU}$ and $(Y_J^\circ)^{\cU}$, with the least
column omitted on each side.
\end{lemma}
\begin{proof}
The bounds $\h\le\Nrm\le\sqrt3\,\h$ identify the underlying
Hilbert and Banach ultrapowers as vector spaces. In each Hilbert
ultrapower, let $A_\Lambda$ be the closed span of the indexed
vectors and let $D_\Lambda=A_\Lambda^\perp$, for $\Lambda=I,J$.
For the ultrapowers with the least column omitted, an indexed
column $F_\beta$ with $\beta>[\min\Lambda]_{\cU}$ can be represented
by a sequence using only columns strictly above $\min\Lambda$.
In either case,
\begin{equation}\label{eq:qultradecomp}
 q_\Lambda^{\cU}(x)=\sup_{\alpha<\beta\text{ in }\Lambda^{\cU}}
 \left|\langle x,E_\alpha\rangle_{\h}+
            \langle x,F_\beta\rangle_{\h}\right|.
\end{equation}
Write $x=[x_s]_{\cU}$, so $q_\Lambda^{\cU}(x)=\lim_{\cU}q_\Lambda(x_s)$.
For each $\alpha=[i_s]_{\cU}<\beta=[j_s]_{\cU}$, the absolute
coordinate sum at $(i_s,j_s)$ is at most $q_\Lambda(x_s)$ on the set in
$\cU$ where $i_s<j_s$. Taking ultralimits and then the supremum
shows that the right-hand side of \eqref{eq:qultradecomp} is at most
$q_\Lambda^{\cU}(x)$. For the reverse inequality, fix $\varepsilon>0$ and choose in each factor
an ordered pair whose absolute coordinate sum exceeds $q_\Lambda(x_s)-\varepsilon$.
The classes of these pairs give a term of size at least
$q_\Lambda^{\cU}(x)-\varepsilon$ in the supremum. Letting $\varepsilon$ tend to zero proves
\eqref{eq:qultradecomp}.

Relabeling by $\sigma$ gives a surjective Hilbert isometry
$A_I\to A_J$ and preserves the displayed supremum. To obtain
an isometry on the complements, choose a bijection of $I$ onto
$J$ that sends the least point to the least point. The associated
Hilbert coordinate isometry has an ultrapower that maps $A_I$
onto $A_J$, so its restriction maps $D_I$ onto $D_J$. Take the
direct sum of this restriction with the relabeling by $\sigma$
on $A_I$. It preserves $\h$ and $q$, hence also the Banach norm.
\end{proof}

\begin{lemma}[Elementary equivalence and the initial successor chain]
\label{lem:successor-inclusion}
For any countable $I,J\equiv(\mathbb N,<)$ one has
$Y_I\equiv Y_J$.  If $c_1,c_2,\ldots$ is the initial successor
chain of $J$, the coordinate map
\[
 Y_{\mathbb N}\longrightarrow Y_J,\qquad
 E_n\mapsto E_{c_n},\quad F_n\mapsto F_{c_n},
\]
is elementary.  These assertions hold in the pure Banach language.
\end{lemma}
\begin{proof}
By the Keisler--Shelah theorem \cite{Shelah1971}, there is a common
ultrafilter $\cU$ and an order isomorphism $\sigma:I^{\cU}\to J^{\cU}$.
Lemma~\ref{lem:ultra-relabel} gives a Banach ultrapower isometry.
By \L o\'s's theorem, $Y_I\equiv Y_J$.

Apply this construction with $I=\mathbb N$. The order isomorphism
$\sigma$ preserves the least point and the immediate successor
relation. The diagonal classes $[n]_{\cU}$ and $[c_n]_{\cU}$
are the $(n-1)$st successors of the least points of their
respective ultrapowers, so $\sigma([n]_{\cU})=[c_n]_{\cU}$.
Let $\Psi:Y_{\mathbb N}^{\cU}\to Y_J^{\cU}$ be the Banach ultrapower
isometry, let $d_{\mathbb N},d_J$ be the diagonal embeddings, and
write $\iota$ for the coordinate inclusion. By the action of $\Psi$
on the diagonal coordinates, the maps $\Psi\circ d_{\mathbb N}$ and
$d_J\circ\iota$ agree on $E_n,F_n$ for each $n$. These vectors have
dense linear span in $Y_{\mathbb N}$, so
\[
 \Psi\circ d_{\mathbb N}=d_J\circ\iota.
\]
Since $\Psi$ is an isomorphism in the pure Banach language and both
diagonal embeddings are elementary, $\iota$ is elementary.
\end{proof}

\begin{lemma}[Adding a Hilbert summand]\label{lem:add-hilbert-summand}
For a countable $I\equiv(\mathbb N,<)$ and separable Hilbert spaces
$K$ and $K_1$, the canonical inclusion
$Y_I^\circ\oplus_2 K\to Y_I^\circ\oplus_2(K\oplus_2 K_1)$
is elementary.
\end{lemma}
\begin{proof}
Fix a nonprincipal ultrafilter $\cU$ on $\mathbb N$. In the Hilbert
space underlying $(Y_I^\circ)^{\cU}$, let $A$ be the closed span of
the indexed row vectors and the indexed column vectors above the
least point of $I^{\cU}$, and put $D=A^\perp$. By
\eqref{eq:qultradecomp}, the seminorm depends only on the projection
onto $A$. The diagonal copy of $Y_I^\circ$ lies in $A$.

To embed $(\ell_2(\mathbb N))^{\cU}$ into $D$, use normalized
averages of disjoint blocks of row vectors. Enumerate $I$ as
$i_1,i_2,\dots$ without repetition.
For $n\ge1$ define a Hilbert isometry $V_n$ from
$\ell_2(\mathbb N)$ into the Hilbert space underlying $Y_I^\circ$ by
\[
 V_n e_m=\frac{1}{\sqrt n}\sum_{r=1}^n E_{i_{(m-1)n+r}}
 \qquad(m\ge1).
\]
For fixed $n$ the supports in this display are disjoint and
the vectors have Hilbert norm one. The map
$[x_n]_{\cU}\mapsto[V_nx_n]_{\cU}$ is a Hilbert isometry on
$(\ell_2(\mathbb N))^{\cU}$ because each $V_n$ is a Hilbert isometry.
If $\|x_n\|_2\le C$ for every $n$, then for any indexed
row vector $E_\alpha$ with $\alpha=[j_n]_{\cU}$,
\[
 |\langle V_nx_n,E_{j_n}\rangle_{\h}|\le C/\sqrt n .
\]
The ultralimit is zero, and every column inner product is zero
already in each factor. Hence the image lies in $D$.

As $K_1$ is separable, $K_1^{\cU}$ embeds isometrically into
$(\ell_2(\mathbb N))^{\cU}$ and hence into $D$. The Hilbert space
$D$ is infinite-dimensional, so there is a surjective Hilbert
isometry $D\to D\oplus_2 K_1^{\cU}$. Take its direct sum with the
identities on $A$ and $K^{\cU}$, placing the $K_1^{\cU}$ component
last. The resulting map is a Hilbert isometry of the two
ultrapowers. It preserves the seminorm because it fixes $A$, and
it extends the canonical inclusion on the diagonal copies because
it fixes $A$ and $K^{\cU}$. Elementarity of the diagonal maps now
gives the result.
\end{proof}

\begin{theorem}[Classification of separable models]
\label{thm:separable-classification}
The separable models of $\operatorname{Th}(\Y)$ in the pure Banach
language are, up to surjective linear isometry, exactly
\begin{equation}\label{eq:separable-classification}
 Y_I^\circ\oplus_2K,
\end{equation}
where $I\models\operatorname{Th}(\mathbb N,<)$ is a countable
linear order and $K$ is a separable real Hilbert space, possibly
zero.  The Banach isometry class determines the order-isomorphism
type of $I$ and the Hilbert dimension of $K$.  For every such
$I,K$, the inclusion $Y_I^\circ\preccurlyeq
Y_I^\circ\oplus_2K$ is elementary.
\end{theorem}
\begin{proof}
\textit{Step 1. Signed pair vectors and coordinate vectors.}
Let $M$ be a separable model of $\operatorname{Th}(\Y)$.
By the metric Keisler--Shelah theorem \cite[Theorem~5.7]{BBHU08},
there are an ultrafilter $\cU$ on an index set $\Omega$ and an
isomorphism $\theta:M^{\cU}\to \Y^{\cU}$.
Let $d_M:M\to M^{\cU}$ be the diagonal embedding.
The composition $\theta\circ d_M$ is elementary, and we identify $M$
with its image in $Z=\Y^{\cU}$.
On $Z$, the predicates $\h$ and $q$ from Lemma~\ref{lem:intrinsic-h}
are given by ultralimits of their values on $\Y$. Their restrictions
to $M$ are the same predicates computed in $M$.

Put $\Lambda=(\mathbb N,<)^\Omega/\cU$, write $c_1,c_2$ for its
first two points, and let
\[
 \mathcal E_\Lambda=\left\{\pm\frac{E_\alpha+F_\beta}{\sqrt2}:
                       \alpha<\beta\text{ in }\Lambda\right\}
      \subseteq B_2^Z .
\]
Here $E_\alpha=[E_{i_\xi}]_{\cU}$ for $\alpha=[i_\xi]_{\cU}$
with $\xi\in\Omega$, and similarly for $F_\beta$. The vectors
$E_\alpha,F_\beta$ form an orthonormal family because their
inner products are ultralimits of coordinate inner products.
Each member of $\mathcal E_\Lambda$ has Hilbert norm one and
thus Banach norm at most $\sqrt3$. Distinct members have inner
product at most $1/2$ and hence Hilbert distance at least one,
so $\mathcal E_\Lambda$ is closed.
By \eqref{eq:qultradecomp}, for $x\in B_2^Z$ we have
\[
 \inf_{e\in \mathcal E_\Lambda}\h(x-e)^2
 =\h(x)^2+1-2\sup_{e\in \mathcal E_\Lambda}\langle x,e\rangle_{\h}
 =\h(x)^2+1-\sqrt2\,q(x).
\]
Thus the predicate $d_{\mathcal E}$ in \eqref{eq:edge-distance} is the
Hilbert distance to $\mathcal E_\Lambda$. Its zero set
$\mathcal E(Z)$ equals $\mathcal E_\Lambda$.

In $G(Z)=\mathcal E(Z)/\{\pm1\}$, the relation $\parallel$ holds exactly when
the pairs share a row or a column, by the inner product calculation
in Theorem~\ref{thm:order-interpretation}. By Lemma~\ref{lem:graph-indices},
the first row consists of the classes $p_j=(c_1,j)$ for $j>c_1$,
ordered by their column indices.
Theorem~\ref{thm:order-interpretation} gives $G(M)\preccurlyeq G(Z)$
and $P(M)\preccurlyeq P(Z)$.
Put
\[
 J=\{c_1\}\cup\{j>c_1:p_j\in P(M)\}\subseteq \Lambda.
\]

The least point $p_{c_2}$ of $P(Z)$ belongs to $P(M)$, so
$c_2\in J$. For $j\in J\setminus\{c_1\}$, the successor
of $p_j$ in $P(Z)$ is $p_{s(j)}$, where $s$ is the successor
in $\Lambda$. This point belongs to $P(M)$ by elementarity.
Since $s(c_1)=c_2$, it follows that $s(J)\subseteq
J\setminus\{c_1\}$. The predecessor of $c_2$ in $\Lambda$
is $c_1$. For each $j\in J$ with $j>c_2$, the predecessor
of $p_j$ in $P(Z)$ belongs to $P(M)$ by elementarity.
Its column is the predecessor of $j$ in $\Lambda$, so that
predecessor belongs to $J$. Thus $s$ is an order isomorphism
from $J$ onto $J\setminus\{c_1\}\cong P(M)$, and
$J\equiv(\mathbb N,<)$.

For every $j\in J\setminus\{c_1\}$ one of the two signed
representatives of $p_j$ belongs to $M$. Closure under
scalar multiplication gives $E_{c_1}+F_j\in M$. The average
of $n$ such vectors for distinct $j$ differs from $E_{c_1}$
by the average of $n$ orthonormal vectors, of Hilbert norm
$1/\sqrt n$. Since $M$ is closed for the equivalent Hilbert
norm, $E_{c_1}$ belongs to $M$, and subtraction gives
$F_j\in M$ for all $j\in J\setminus\{c_1\}$.
This orthonormal family lies in a separable Hilbert space,
so $J$ is countable.

\emph{A cell $(i,j)$ of $G(Z)$ belongs to $G(M)$ exactly when
$i,j\in J$.} For a cell $e$ outside the first row, let $A(e)$ be
the set of $p\in P$ for which there is a cell $g\notin P$ with
$g\parallel e$ and $g\parallel p$, as in Lemma~\ref{lem:graph-indices}.
If $e=(i,j)$ belongs to $G(M)$ and is outside $P(M)$, the lemma
identifies $p_j$ in $G(Z)$ as the unique $p\in P$ with $e\parallel p$
and $p_{s(i)}$ as the least element of $A(e)$. Thus both lie in
$P(M)$ by elementarity and uniqueness in $G(Z)$. Hence $j$ and
$s(i)$ belong to $J$, and so does $i$ by predecessor closure.
If a cell lies in $P(M)$, its row is $c_1$ and its column belongs to
$J$ by definition.

Conversely, fix $i<j$ in $J$. If $i=c_1$, then $(i,j)=p_j\in P(M)$.
If $i>c_1$, both $p_j$ and $p_{s(i)}$ belong to $P(M)$.
In $G(Z)$ the cell $(i,j)$ witnesses the formula
\[
 \exists e\notin P\,\bigl(e\parallel p_j\wedge p_{s(i)}\in A(e)\wedge
          \forall p'\in A(e)\ \neg(p'\prec p_{s(i)})\bigr).
\]
By elementarity there is a witness $(i',j')$ in $G(M)$.
Since $G(M)\preccurlyeq G(Z)$, the conjunction holds for the same
witness in $G(Z)$, where the universal quantifier ranges over all
of $A((i',j'))$ in $G(Z)$. By Lemma~\ref{lem:graph-indices},
$(i',j')\parallel p_j$ gives $j'=j$. The least element of $A((i',j'))$
is both $p_{s(i')}$ and $p_{s(i)}$, so $s(i')=s(i)$ and $i'=i$.
Thus $(i,j)\in G(M)$ in both cases. Each such class has a
representative in $\mathcal E(M)$, and its negative is the other
representative. Hence
\[
 \mathcal E(M)=\left\{\pm\frac{E_i+F_j}{\sqrt2}:
                                    i<j\text{ in }J\right\}.
\]
In particular, for each $i>c_1$ in $J$ the cell $(i,s(i))$ gives
$E_i+F_{s(i)}\in M$. Subtracting $F_{s(i)}$ gives $E_i\in M$.

\textit{Step 2. Orthogonal decomposition.}
Let $A$ be the closed Hilbert span of the vectors $E_i$ for
$i\in J$ and $F_j$ for $j\in J\setminus\{c_1\}$, and put
$K=M\cap A^\perp$. For $x\in B_2^M$, Lemma~\ref{lem:edge-set}
gives $d_{\mathcal E}(x)=\dist_{\h}(x,\mathcal E(M))$. Expanding the squared
distance to these Hilbert unit vectors and using both signs gives
\[
 d_{\mathcal E}(x)^2=\h(x)^2+1-2\sup_{e\in\mathcal E(M)}
                              |\langle x,e\rangle_{\h}|.
\]
Comparing with \eqref{eq:edge-distance} yields
\[
 q(x)=\sqrt2\sup_{e\in\mathcal E(M)}|\langle x,e\rangle_{\h}|,
\]
which holds on all of $M$ by homogeneity. Every signed pair vector
lies in $A$, so $q(a+k)=q(a)$ for $a\in A$ and $k\in K$.
The coordinate Hilbert isometry from the underlying space
of $Y_J^\circ$ onto $A$ identifies $q_J$ with $q|_A$ by the
description of $\mathcal E(M)$. Hence
\[
 \Nrm(a+k)^2=\h(a)^2+q(a)^2+\h(k)^2.
\]
The Hilbert decomposition of $M$ into $A$ and $K$ therefore gives
a surjective linear isometry $Y_J^\circ\oplus_2 K\to M$, with
$K$ carrying the Hilbert norm $\h|_K$.

The kernel of $q_J$ on $Y_J^\circ$ is zero.  Indeed, if
$u_i+v_j=0$ for all $i<j$, choose $j$ above any given $i,k$
to get $u_i=u_k$. Square summability forces this common value
to be zero. For each $j\in J$ above its least point $c_1$, the pair
$(c_1,j)$ gives $v_j=0$. Hence $K=\ker q$ and $\Nrm|_K=\h|_K$.

\textit{Step 3. Elementary inclusion after adding a Hilbert summand.}
By Lemma~\ref{lem:add-hilbert-summand} with the initial Hilbert
summand zero, $Y_I^\circ\preccurlyeq Y_I^\circ\oplus_2 K$.
Taking $K=\mathbb R$ gives $Y_I^\circ\preccurlyeq Y_I$ under
the decomposition along $F_{c_I}$. Since $Y_I\equiv\Y$ by
Lemma~\ref{lem:successor-inclusion}, every
$Y_I^\circ\oplus_2K$ is a model of $\operatorname{Th}(\Y)$.

\textit{Step 4. Isometry invariants.}
By Lemma~\ref{lem:intrinsic-h}, on $Y_I^\circ\oplus_2K$ the intrinsic
$\h$ is the Hilbert sum norm and $q(u,v,k)=q_I(u,v)$.
For $x=(u,v,k)$, the inner products with the signed pair vectors satisfy
\[
 \sup_{\substack{i<j\\ \varepsilon\in\{-1,1\}}}
 \left\langle x,\frac{\varepsilon(E_i+F_j)}{\sqrt2}\right\rangle_{\h}
 =\frac{1}{\sqrt2}\sup_{i<j}|u_i+v_j|
 =\frac{q_I(u,v)}{\sqrt2}.
\]
These vectors have Hilbert norm one and pairwise Hilbert distance at
least one, so they form a closed set. Expanding the squared distance
to them and using the supremum above gives \eqref{eq:edge-distance}
on $B_2$. Thus these vectors form the intrinsic set $\mathcal E$
in this model. By Lemma~\ref{lem:graph-indices}, the
first-row order on their sign classes is $I\setminus\{c_I\}$, which
is order-isomorphic to $I$ via successor. Also, $K=\ker q$ by the
kernel computation in Step~2. A surjective linear Banach isometry
preserves the interpreted order and the Hilbert norm on $\ker q$,
so it preserves the order-isomorphism type of $I$ and the Hilbert
dimension of $K$. Conversely, an order isomorphism $I\to J$
together with a Hilbert isometry
$K\to K'$, gives a surjective Banach isometry between the
corresponding sums.
\end{proof}

For a countable model $I$ of $\operatorname{Th}(\mathbb N,<)$, define
\[
 i\sim j\quad\Longleftrightarrow\quad
 \text{the closed interval with endpoints $i$ and $j$ is finite}.
\]
\begin{lemma}[Decomposition into convex blocks]\label{lem:index-block-form}
The relation $\sim$ is an equivalence relation with convex classes.
The initial class has order type $\omega$ and every other class
has order type $\mathbb Z$. If $B_I$ is the order of the noninitial
classes, then
\begin{equation}\label{eq:index-block-form}
 I\cong\omega+\sum_{b\in B_I}\mathbb Z.
\end{equation}
\end{lemma}
\begin{proof}
Finite intervals are closed under taking subintervals and under
joining two intervals with a common endpoint. Thus $\sim$ is an
equivalence relation and its classes are convex. Two points lie
in the same class exactly when one is obtained from the other
by finitely many successor or predecessor steps.

The class of the least point consists of that point and all its
finite successors, so it has order type $\omega$. Any point
outside this class has predecessors of every finite depth,
since otherwise finitely many predecessor steps would reach
the least point. It also has successors of every finite depth.
Its class therefore has order type $\mathbb Z$. The order sum
in \eqref{eq:index-block-form} lists the classes in their order
in $I$.
\end{proof}

\begin{corollary}[Separable spectrum]\label{cor:separable-spectrum}
There are exactly $2^{\aleph_0}$ pairwise nonisometric separable
models of $\operatorname{Th}(\Y)$.
\end{corollary}
\begin{proof}
For each countable order $B$, put
$I_B=\omega+\sum_{b\in B}\mathbb Z$ in lexicographic order.
The order $I_B$ has a least element and no greatest element,
every element has an immediate successor, and every nonleast
element has an immediate predecessor. Thus
$I_B\equiv(\mathbb N,<)$ by Lemma~\ref{lem:ord}(\ref{item:ord-equivalence}).
The order $B_{I_B}$ of its noninitial blocks is isomorphic to $B$.
To obtain $2^{\aleph_0}$ distinct
choices of $B$ up to isomorphism, let $A\subseteq\mathbb N$ and
put $B_A=\sum_{n\in\mathbb N}(C_n+\mathbb Q)$, where $C_n$ is a
finite chain of size $2$ if $n\in A$ and size $3$ otherwise.
The nonsingleton classes for the finite-interval relation in
$B_A$ are precisely the chains $C_n$, since the interval between
two distinct points that do not lie in the same $C_n$ contains
infinitely many points from one of the copies of $\mathbb Q$.
The chains occur in order type
$\omega$, so their sizes recover $A$ from the order type of $B_A$.
The corresponding $Y_{I_{B_A}}^\circ$ are pairwise nonisometric by
Theorem~\ref{thm:separable-classification}. For the upper bound,
there are at most $2^{\aleph_0}$ countable orders up to isomorphism,
since an order on $\mathbb N$ is a subset of $\mathbb N^2$.
There are only countably many possible dimensions of $K$.
\end{proof}

\begin{corollary}[Prime model]\label{cor:prime}
$Y_{\mathbb N}^\circ$ is the prime model of
$\operatorname{Th}(\Y)$.
\end{corollary}
\begin{proof}
For an arbitrary $M\models\operatorname{Th}(\Y)$,
choose a separable elementary submodel $M_0\preccurlyeq M$ by
downward L\"owenheim--Skolem \cite[Proposition~7.3]{BBHU08}
in the countable Banach language.
Theorem~\ref{thm:separable-classification} identifies it with
$Y_I^\circ\oplus_2K$. The successor-chain coordinate
inclusion $Y_{\mathbb N}\to Y_I$ is elementary by
Lemma~\ref{lem:successor-inclusion}. The inclusion sends $F_1$
to $F_{c_I}$ and sends every other column vector to a column
with index greater than $c_I$. Thus the image of a vector
with zero $F_1$-coordinate has zero $F_{c_I}$-coordinate, and
$Y_{\mathbb N}^\circ$ maps into $Y_I^\circ$. Since
$Y_{\mathbb N}^\circ\preccurlyeq Y_{\mathbb N}$ and
$Y_I^\circ\preccurlyeq Y_I$, restriction gives an elementary
map $Y_{\mathbb N}^\circ\to Y_I^\circ$. Compose with
$Y_I^\circ\preccurlyeq M_0\preccurlyeq M$.

Any other prime model embeds elementarily into
$Y_{\mathbb N}^\circ$ and is therefore separable. By
Lemma~\ref{lem:intrinsic-h} its kernel of $q$ embeds into
the zero kernel of $Y_{\mathbb N}^\circ$. By
Theorem~\ref{thm:order-interpretation} its interpreted order
embeds elementarily into $(\mathbb N,<)$. The image contains
the least element and is closed under successor, so this order
is isomorphic to $(\mathbb N,<)$. The classification in
Theorem~\ref{thm:separable-classification} gives an isometry
with $Y_{\mathbb N}^\circ$.
\end{proof}

The space $Y$ itself is not prime. Its kernel of $q$ is
$\mathbb R F_1$, so preservation of $q$ rules out an elementary
embedding into $Y_{\mathbb N}^\circ$, where the kernel is zero.

\subsection{Elementary embeddings and automorphism groups}
\label{sec:elementary-maps}
For an order embedding $f:I\to J$ sending the least point to the
least point, define the coordinate map
$f_*:Y_I^\circ\to Y_J^\circ$ by
$E_i\mapsto E_{f(i)}$ and
$F_j\mapsto F_{f(j)}$ for $j>c_I$, extending by linearity
and continuity in the Hilbert norms.

\begin{theorem}[Elementary embeddings of separable models]
\label{thm:embedding-normal-form}
Every elementary embedding
\[
 \Phi:Y_I^\circ\oplus_2K\longrightarrow Y_J^\circ\oplus_2K'
\]
has a unique presentation
\begin{equation}\label{eq:embedding-normal}
 \Phi(z,k)=(\varepsilon f_*z,Vk),
\end{equation}
where $f:I\to J$ is an elementary order embedding,
$\varepsilon\in\{-1,1\}$, and $V:K\to K'$ is a linear
Hilbert isometric embedding.  Conversely each such triple defines
an elementary Banach embedding.
\end{theorem}
\begin{proof}
An elementary Banach embedding preserves addition, rational
scalar multiplication, and the norm, so it is a linear isometry
over the reals by continuity. It also preserves the intrinsic Hilbert norm and $q$ by
Lemma~\ref{lem:intrinsic-h}. It maps signed pair vectors to signed pair vectors
by Lemma~\ref{lem:edge-set} and induces an elementary embedding of
the interpreted first-row orders by Theorem~\ref{thm:order-interpretation}.
Write $g(j)$ for the column of the image of
the first-row class $(c_I,j)$. By
Theorem~\ref{thm:order-interpretation}, this defines an elementary embedding
$g:I\setminus\{c_I\}\to J\setminus\{c_J\}$.

Write $s_I,s_J$ for the successor maps on the full orders.
Each is increasing and maps its order bijectively onto the
tail, because every nonleast point has an immediate predecessor.
Hence
\[
 f=s_J^{-1}\circ g\circ s_I:I\longrightarrow J
\]
is elementary. The map $g$ preserves the least points of the
tails, so $g(s_I(c_I))=s_J(c_J)$ and $f(c_I)=c_J$. For $i>c_I$,
the successor in the tail is the same as in the full order. As
$g$ preserves this successor, $g(s_I(i))=s_J(g(i))$, giving
$f(i)=g(i)$.

Thus the first-row cell $(c_I,j)$ maps to $(c_J,f(j))$.
For a cell $e=(i,j)$ outside the first row, Lemma~\ref{lem:graph-indices}
identifies $p_j$ as the unique $p\in P$ with $e\parallel p$ and
$p_{s_I(i)}$ as the least element of $A(e)$. Their images have
column indices $g(j)$ and $g(s_I(i))$. The image cell therefore has column
$f(j)$ and row $s_J^{-1}(g(s_I(i)))=f(i)$.

On signed representatives this says
$\Phi(e_{ij})=\varepsilon_{ij}e_{f(i)f(j)}$ for signs
$\varepsilon_{ij}\in\{\pm1\}$. If distinct cells $(i,j)$ and $(k,\ell)$
share a row or a column, their images do as well. Preservation
of the Hilbert inner product then gives
\[
 \frac12=\langle e_{ij},e_{k\ell}\rangle_{\h}
 =\langle\Phi(e_{ij}),\Phi(e_{k\ell})\rangle_{\h}
 =\frac{\varepsilon_{ij}\varepsilon_{k\ell}}2.
\]
Thus $\varepsilon_{ij}=\varepsilon_{k\ell}$. The first-row cells share
a row, and every other cell shares a column with a first-row cell, so
there is a single global sign $\varepsilon$.
As $j$ ranges over any sequence of distinct indices above $c_I$,
the averages of $E_{c_I}+F_j=\sqrt2\,e_{c_Ij}$ converge in Hilbert
norm to $E_{c_I}$, since the $F_j$ are orthonormal. The closed span
of the pair vectors therefore contains $E_{c_I}$ and, by subtraction,
each $F_j$ with $j>c_I$. For each $i$ choose $j>i$ and subtract $F_j$
from $\sqrt2\,e_{ij}$ to obtain $E_i$. Thus the pair vectors span
$Y_I^\circ$ densely, and
$\Phi|_{Y_I^\circ}=\varepsilon f_*$. Since $\Phi$ preserves $q$,
it maps $K=\ker q$ into $K'$, and its restriction $V$ is a
Hilbert isometric embedding. The maps $f$, $V$, and the sign
$\varepsilon$ are determined by $\Phi$.

For the converse, first construct the map
$(z,k)\mapsto(f_*z,k)$ into $Y_J^\circ\oplus_2 K$.
Expand the two index orders by constants naming every point of
$I$, interpreted as $i$ in $I$ and $f(i)$ in $J$.  Since $f$ is
an \emph{elementary} order embedding, the two expanded orders
are elementarily equivalent.  A common-ultrafilter
Keisler--Shelah isomorphism therefore gives
\[
 \sigma:I^{\cU}\stackrel{\cong}{\longrightarrow}J^{\cU},
 \qquad\sigma([i]_{\cU})=[f(i)]_{\cU}
       \quad\text{for every }i\in I.
\]
By Lemma~\ref{lem:ultra-relabel}, $\sigma$ induces a surjective linear
isometry $\widehat\sigma:(Y_I^\circ)^{\cU}\to(Y_J^\circ)^{\cU}$.
It sends the diagonal $E_i$ to $E_{f(i)}$ for every $i\in I$ and
the diagonal $F_j$ to $F_{f(j)}$ for every $j>c_I$.
The source coordinate vectors span the diagonal copy of $Y_I^\circ$
densely, so $\widehat\sigma$ acts there by $f_*$.

For $L=I,J$, taking ultralimits in the direct-sum norm identifies
$(Y_L^\circ\oplus_2 K)^{\cU}$ with
$(Y_L^\circ)^{\cU}\oplus_2 K^{\cU}$ by
$[(z_\xi,k_\xi)]_{\cU}\mapsto([z_\xi]_{\cU},[k_\xi]_{\cU})$.
Write $d_L$ for the diagonal embedding into this ultrapower. Then
\[
 (\widehat\sigma\oplus\mathrm{id}_{K^{\cU}})\,d_I(z,k)
 =d_J(f_*z,k).
\]
Both diagonal embeddings are elementary. The surjective linear isometry
$\widehat\sigma\oplus\mathrm{id}_{K^{\cU}}$ is an isomorphism in the
pure Banach language. The identity therefore gives the elementary embedding
\[
 (z,k)\longmapsto (f_*z,k):
   Y_I^\circ\oplus_2 K\preccurlyeq Y_J^\circ\oplus_2 K.
\]
Next consider the map $(z,k)\mapsto(z,Vk)$ from $Y_J^\circ\oplus_2K$
to $Y_J^\circ\oplus_2K'$. Identify $K$ with $V(K)\subseteq K'$.
The orthogonal complement of $V(K)$ in $K'$ is separable, so
Lemma~\ref{lem:add-hilbert-summand} gives elementarity of this map.
The two
elementary embeddings compose as
\[
 (z,k)\longmapsto(f_*z,k)\longmapsto(f_*z,Vk).
\]
Multiplying the coordinate summand of the target by
$\varepsilon$ preserves $\h$ and $q$ and gives a Banach
automorphism. Composing with it yields \eqref{eq:embedding-normal}.
\end{proof}

For the renorming of $c_0$ due to Ferenczi and Rosendal, each
surjective linear isometry similarly acts on the coordinate
basis by a graph automorphism and a common sign
\cite[Section~3.2, Lemma~11]{FerencziRosendal}.

\begin{corollary}[Algebraic and definable closure]\label{cor:algebraic-closure}
Let $Y_0$ be the image of an elementary embedding of $Y_{\mathbb N}^\circ$
into a monster model $\mathfrak U$ of $\operatorname{Th}(\Y)$. In the
Banach sorts,
\[
 \operatorname{acl}(\varnothing)=Y_0,
 \qquad \operatorname{dcl}(\varnothing)=\{0\}.
\]
For every nonzero $a\in Y_0$, one has
$\operatorname{dcl}(a)=\operatorname{acl}(a)=Y_0$.
\end{corollary}
\begin{proof}
Fix $\sigma\in\aut(\mathfrak U)$. By downward L\"owenheim--Skolem,
countable dense subsets of $Y_0$ and $\sigma(Y_0)$ lie in a common
separable elementary submodel $M$ of $\mathfrak U$. Since $M$ is
closed, it contains both images. Since $M\preccurlyeq\mathfrak U$,
the inclusion of $Y_0$ and the restriction of $\sigma$ to $Y_0$ are
elementary embeddings into $M$. Identify $M$ with $Y_I^\circ\oplus_2 K$
by Theorem~\ref{thm:separable-classification}. The order embeddings
associated to the two Banach embeddings by
Theorem~\ref{thm:embedding-normal-form} agree, since every elementary
embedding of $(\mathbb N,<)$ into $I$ takes $1$ to the least point and
preserves immediate successors. Thus $\sigma|_{Y_0}$ is multiplication
by a common sign. Global negation is a Banach automorphism, so the
orbit of each $x\in Y_0$ is $\{x,-x\}$.

These orbits are compact, so $Y_0\subseteq\operatorname{acl}(\varnothing)$
\cite[Exercise~10.8]{BBHU08}. Since $Y_0\preccurlyeq\mathfrak U$,
every vector algebraic over the empty set belongs to $Y_0$
\cite[Corollary~10.5]{BBHU08}. Zero is the only
fixed point of global negation, and is named in the Banach language.
Thus $\operatorname{dcl}(\varnothing)=\{0\}$. If $a\in Y_0$ is nonzero,
every automorphism fixing $a$ has positive sign on $Y_0$. It follows
that $Y_0\subseteq\operatorname{dcl}(a)$. Since $Y_0$ is an elementary
submodel containing $a$, we have
$\operatorname{dcl}(a)\subseteq\operatorname{acl}(a)\subseteq Y_0$
\cite[Corollary~10.5 and Exercise~10.7]{BBHU08}.
\end{proof}

\begin{corollary}[Automorphism group]\label{cor:automorphism-group}
For each separable model $Y_I^\circ\oplus_2K$ there is an
isomorphism of topological groups
\begin{equation}\label{eq:aut-group}
 \operatorname{Aut}(Y_I^\circ\oplus_2K)
 \cong\{\pm1\}\times\operatorname{Aut}(I,<)\times O(K).
\end{equation}
The two automorphism groups have the topologies of pointwise
convergence in the Banach norm and in the discrete order, respectively.
The group
$O(K)$ of surjective linear isometries of $K$ carries the strong
operator topology.
\end{corollary}
\begin{proof}
A map of the form \eqref{eq:embedding-normal} is onto precisely
when $f$ and $V$ are onto.  Composition acts componentwise on
the signs, orders, and Hilbert maps. Signed row vectors at
distinct indices have Hilbert distance $\sqrt2$, and opposite
signs at the same index have Hilbert distance $2$. Since
$\Nrm\ge\h$, pointwise norm convergence on $E_i$ forces
eventual agreement of the sign and of $f(i)$ for every $i$.
Evaluation on $K$ gives strong operator
convergence. Conversely, eventual agreement on signs and on
the images of finitely many indices gives norm convergence on every
finitely supported coordinate vector. Density and the norm-one
bound on the automorphisms give convergence on $Y_I^\circ$.
Strong convergence on $K$ then gives convergence on the sum.
\end{proof}

\begin{corollary}[Automorphisms of the index order]\label{cor:order-blocks}
Let $B_I$ be the order of the $\mathbb Z$-blocks in \eqref{eq:index-block-form}.
After choosing an origin in each $\mathbb Z$-block, there is
an isomorphism of topological groups
\begin{equation}\label{eq:block-automorphisms}
 \operatorname{Aut}(I,<)
   \cong\mathbb Z^{B_I}\rtimes\operatorname{Aut}(B_I,<),
\end{equation}
where the shift group carries the product of discrete topologies
and the order-automorphism group carries pointwise convergence.
\end{corollary}
\begin{proof}
An order automorphism preserves the least point, the successor map,
and the finite-interval equivalence relation. It fixes the initial
$\omega$-block pointwise, induces an automorphism of $B_I$, and maps each
$\mathbb Z$-block to its image by an integer translation.
For $\sigma\in\operatorname{Aut}(B_I,<)$
and $a\in\mathbb Z^{B_I}$, index the shift by the target block
by letting $(a,\sigma)$ act on the noninitial blocks by
\[
 (b,n)\longmapsto(\sigma(b),n+a(\sigma(b))).
\]
Every automorphism of $I$ has this form for a unique
$(a,\sigma)$, and the action gives the group law
\[
 (a,\sigma)(d,\tau)=(a+\sigma\cdot d,\sigma\circ\tau),
 \qquad (\sigma\cdot d)(c)=d(\sigma^{-1}(c)).
\]

Pointwise convergence of the automorphisms of $I$ gives eventual
agreement on the image of each $(b,0)$. This implies eventual
agreement on $\sigma(b)$, and taking $b=\sigma^{-1}(c)$ gives
eventual agreement on the shift $a(c)$ for each target block $c$.
Conversely, eventual agreement on the image of each block and
on the shift at each target block gives eventual agreement at
every $(b,n)$ by the action formula.
\end{proof}

\begin{corollary}[Elementary embeddings of countable index orders]
\label{cor:index-embeddings}
After choosing origins in the $\mathbb Z$-blocks of $I$ and $J$,
every elementary order embedding $f:I\to J$ sends the $n$th
successor of $c_I$ to the $n$th successor of $c_J$ for every
$n\ge0$. On the noninitial blocks it has the form
\[
 (b,n)\longmapsto\bigl(\beta(b),n+a(b)\bigr)
 \qquad(b\in B_I,\ n\in\mathbb Z),
\]
where $\beta:B_I\to B_J$ is an order embedding and
$a:B_I\to\mathbb Z$ is an arbitrary function.  Conversely every
such pair defines an elementary order embedding.
\end{corollary}
\begin{proof}
An elementary embedding preserves least points and immediate
successors and predecessors, so it sends the initial successor
chain onto the initial successor chain. The initial gaps of
points in noninitial blocks are infinite, so their images are
outside the initial block. Since the embedding commutes with
successor and predecessor, it maps each $\mathbb Z$-block onto
a $\mathbb Z$-block by an integer translation.
Points in distinct blocks have infinitely many points between
them. Elementarity preserves each finite lower bound on that
number, so their images lie in distinct blocks. The induced map
$\beta$ is increasing.
Conversely a map of the displayed form preserves the exact
finite gaps between points of the same block and the fact that
gaps between different blocks are larger than every finite
integer. It preserves initial gaps exactly on the initial block,
and these gaps are infinite on every other block. By
Lemma~\ref{lem:ord}(\ref{item:ord-types}), corresponding finite tuples have identical
complete order types.  Thus the map is elementary.
\end{proof}

\subsection{Strong separability of bounded norm types}
\label{sec:strong-types}

Let $B_M$ be the closed unit ball of a real Banach space $M$.
For $y\in B_M$ put $p_y(x)=\|x+y\|_M$ on $B_M$. The space
$S_{\mathrm{norm}}(M)$ of bounded norm types is the pointwise closure
of these functions in $[0,2]^{B_M}$. Its elements are exactly the
functions $x\mapsto\|x+y\|_{M'}$ on $B_M$ with $M'\succeq M$ and
$y\in B_{M'}$. Indeed, elementarity gives approximations on each
finite subset of $B_M$, and compactness realizes the conditions given
by a pointwise limit. Following \cite[Definitions~3.4 and~4.15]{KhanakiBanach}
and Lupini \cite[Section~3]{LupiniNormFormulas2026}, we say that the types
are \emph{strongly separable} if $S_{\mathrm{norm}}(M)$ is separable for
the metric $\rho(p,p')=\sup_{x\in B_M}|p(x)-p'(x)|$.

The calculation of the types applies more generally to any countable
linear order $I$ with a least element and no greatest element. Use the
same coordinate norm and omit the least second coordinate as in
Section~\ref{sec:separable-classification}. Write $I_\infty$ for the
set of $j\in I$ with infinitely many points below $j$. Every $i\in I$
has infinitely many points above it.

If $z_n$ converges weakly to zero, each fixed coordinate tends to
zero, but a column of $z_n$ with a varying index can still contribute
to $q_I(s+z_n)$ together with a fixed row of $s$. The possible
contributions from that row depend on the supremum and infimum of the
column coefficients above it. The same holds with rows and columns
interchanged when the column has infinitely many predecessors.

\begin{lemma}[Limits of the ordered-pair seminorm]\label{lem:weak-pair-limit}
Let $I$ be a countable linear order with a least element and no
greatest element. Let $z_n=(a^{(n)},b^{(n)})$ converge weakly to zero
in the real Hilbert space $H_I$. There is a subsequence on which $q_I(z_n)\to \kappa$
and the following limits exist:
\begin{align*}
 b^+(i)&=\lim_n\sup_{j>i}b_j^{(n)},&
 b^-(i)&=\lim_n\inf_{j>i}b_j^{(n)} &&(i\in I),\\
 a^+(j)&=\lim_n\sup_{i<j}a_i^{(n)},&
 a^-(j)&=\lim_n\inf_{i<j}a_i^{(n)} &&(j\in I_\infty).
\end{align*}
The functions $b^+,b^-$ give the limits of the upper and lower
bounds on column coefficients to the right of a row. The functions
$a^+,a^-$ give the corresponding limits for row coefficients to the
left of a column. All four functions are bounded and monotone.
On the same subsequence,
for every $s=(u,v)\in H_I$,
\begin{equation}\label{eq:limit-seminorm}
 \lim_n q_I(s+z_n)=Q(s):=
 \max\left\{q_I(s),\,\kappa,\,
 \sup_{\substack{i\in I\\\sigma\in\{+,-\}}}|u_i+b^\sigma(i)|,\,
 \sup_{\substack{j\in I_\infty\\\sigma\in\{+,-\}}}|v_j+a^\sigma(j)|\right\}.
\end{equation}
The last supremum is omitted when $I_\infty$ is empty.
\end{lemma}
\begin{proof}
Weak convergence bounds the Hilbert norms of the $z_n$. Countability
of $I$ gives a single subsequence on which all the indicated real
sequences converge. As $i$ increases, the set $\{j\in I:j>i\}$
shrinks, so $b^+$ is nonincreasing and $b^-$ is nondecreasing.
As $j$ increases in $I_\infty$, the set $\{i\in I:i<j\}$ grows, so
$a^+$ is nondecreasing and $a^-$ is nonincreasing. On any infinite
set of indices, the coordinate values of a fixed Hilbert vector have
zero in their closure. Thus $b^-(i)\le0\le b^+(i)$ and $a^-(j)\le0\le a^+(j)$.
For fixed $i$ and any $k>i$, we have
$|b_k^{(n)}|\le q_I(z_n)+|a_i^{(n)}|$. The last term tends to
zero, so $|b^\pm(i)|\le \kappa$. Fixing $j$ instead gives $|a^\pm(j)|\le \kappa$.

For a finite $S\subset I$, let $q_{I\setminus S}$ restrict the
supremum to pairs with both endpoints outside $S$. We claim that
\begin{equation}\label{eq:delete-finite-endpoints}
 q_{I\setminus S}(z_n)\longrightarrow \kappa .
\end{equation}
The upper bound follows from $q_{I\setminus S}\le q_I$, which also
settles $\kappa=0$. Suppose $\kappa>0$ and take pairs $(i_n,j_n)$ whose
contributions differ from $q_I(z_n)$ by at most $1/n$. The
coefficients on $S$ tend uniformly to zero, so eventually at most
one endpoint of such a pair belongs to $S$.
Fix $\eta>0$ and a bound $C$ on the Hilbert norms.
If $i_n\in S$, choose an integer $m>|S|+C^2/\eta^2$. At most $m+1$
points of $I$ have at most $m$ predecessors. On this finite set the
second coefficients also tend uniformly to zero. Thus for all
sufficiently large $n$ with $i_n\in S$, the point $j_n$ has more
than $m$ predecessors. At most $C^2/\eta^2$ first coefficients
have modulus greater than $\eta$, so some $i'<j_n$ outside $S$
has $|a_{i'}^{(n)}|\le\eta$. Replacing $i_n$ by $i'$ costs at
most $\eta+|a_{i_n}^{(n)}|$. If $j_n\in S$, the infinitely many
indices above $i_n$ contain a $j'\notin S$ with
$|b_{j'}^{(n)}|\le\eta$, at a cost of at most $\eta+|b_{j_n}^{(n)}|$.
Pairs outside $S$ already contribute to the restricted supremum.
Its lower limit is therefore at least $\kappa-\eta$, proving
\eqref{eq:delete-finite-endpoints}.

First suppose $s$ has finite support. Take a finite $S$ containing
both supports and all predecessors of each $j\in\supp(v)$ with
finitely many predecessors. Zeroing the coordinates of $z_n$ on $S$
changes $q_I(s+z_n)$ by at most $\sqrt2\,\|P_S z_n\|_2\to0$, where
$P_S$ is the coordinate projection in both copies. Consider the
remaining sums according to the positions of their endpoints.
\begin{itemize}
\item With both endpoints in $S$ the contributions are $|u_i+v_j|$.
      With both endpoints outside $S$ their supremum tends to $\kappa$
      by \eqref{eq:delete-finite-endpoints}.
\item For $i\in S$, the supremum over $j>i$ outside $S$ tends to
      $\max_{\sigma}|u_i+b^\sigma(i)|$. Deleting finitely many
      entries tending to zero preserves the limiting extrema over $j>i$,
      since the remaining supremum is at least zero and the remaining
      infimum is at most zero.
\item For $j\in S\cap I_\infty$, the same argument over $i<j$ outside
      $S$ gives $\max_{\sigma}|v_j+a^\sigma(j)|$. If $j\in S\setminus I_\infty$
      and $v_j\ne0$, all its predecessors are in $S$. If $v_j=0$,
      the finitely many terms with $i<j$ outside $S$ tend to zero.
\end{itemize}
Since $S$ is finite, $q_I(s+z_n)$ converges to the maximum $D$ of
these limiting contributions. For each fixed pair $i<j$, weak
convergence of $z_n$ to zero gives
\[
 |u_i+v_j|
 =\lim_n |u_i+a_i^{(n)}+v_j+b_j^{(n)}|
 \le \lim_n q_I(s+z_n)=D.
\]
Taking the supremum over $i<j$ gives $q_I(s)\le D$. Since $q_I(s)$
also bounds every contribution from a pair inside $S$, replacing the
maximum over those pairs by $q_I(s)$ gives the same value $D$.
Adding indices outside the supports to the two suprema in
\eqref{eq:limit-seminorm} contributes at most $\kappa$. This proves
the identity for finitely supported $s$. Both $q_I(s+z_n)$ and
$Q(s)$ are $\sqrt2$-Lipschitz in $s$ for the Hilbert norm, uniformly
in $n$. Density gives the identity for every $s$.
\end{proof}

\begin{proposition}[Representation of bounded norm types]
\label{prop:bounded-type-representation}
Let $M=Y_I^\circ\oplus_2 K$ be a separable model, and let $x_I$
denote the $Y_I^\circ$-component of $x\in M$. For every
$p\in S_{\mathrm{norm}}(M)$ there are $y\in M$, $r_0,\kappa\ge0$, and
monotone functions $\mathbf f=(b^+,b^-,a^+,a^-)$ on the domains
in Lemma~\ref{lem:weak-pair-limit} such that
\begin{equation}\label{eq:bounded-type-representation}
 p(x)=\bigl(\h(x+y)^2+r_0^2+Q(x_I+y_I)^2\bigr)^{1/2}
 \qquad(x\in B_M).
\end{equation}
Here $Q$ is the function of $\kappa$ and $\mathbf f$ given by
\eqref{eq:limit-seminorm}. The parameters satisfy
$\h(y)^2+r_0^2\le1$, $\kappa\le\sqrt2\,r_0$, and $\|f_k\|_\infty\le r_0$
for each component $f_k$ of $\mathbf f$.
\end{proposition}
\begin{proof}
Every $p\in S_{\mathrm{norm}}(M)$ is a pointwise limit of a
sequence $p_{y_n}$ with $y_n\in B_M$.
Indeed, successive approximations on finite subsets of a countable
dense set suffice, since the functions and their limits are
$1$-Lipschitz. Passing to a subsequence, $y_n$ converges weakly
in the underlying Hilbert space to some $y$, and $\h(y_n-y)\to r_0$.
For each $x$, expansion of the squared Hilbert norm gives
$\h(x+y_n)^2\to\h(x+y)^2+r_0^2$. Apply Lemma~\ref{lem:weak-pair-limit}
to the sequence $(y_n-y)_I$. The resulting function $Q$ gives
\eqref{eq:bounded-type-representation}.
At $x=0$ the Hilbert identity gives $\h(y)^2+r_0^2\le1$.
Since $q_I\le\sqrt2\,\h$, we have $\kappa\le\sqrt2\,r_0$. Each coefficient
of $(y_n-y)_I$ is bounded in magnitude by $\h(y_n-y)$, so each
component of $\mathbf f$ is bounded by $r_0\le1$ in absolute value.
\end{proof}

\begin{theorem}[Strong separability and the order of the blocks]
\label{thm:strong-types}
Let $M=Y_I^\circ\oplus_2 K$ be a separable model of $\operatorname{Th}(Y)$.
\begin{enumerate}[label=(\alph*),ref=\alph*,leftmargin=3em,beginpenalty=10000]
\item\label{item:strong-types-criterion}
The bounded norm types over $M$ are strongly separable if and only if
the order $B_I$ of its noninitial $\mathbb Z$-blocks is scattered,
that is, contains no copy of $(\mathbb Q,<)$.
\item\label{item:strong-types-orders}
Scatteredness of $B_I$ is equivalent to scatteredness of $I$
and to $B_I$ having only countably many initial segments.
\item\label{item:strong-types-separation}
If $B_I$ contains a copy of $(\mathbb Q,<)$, there are continuum
many bounded norm types at pairwise distance at least
$\sqrt3-\sqrt{3/2}$.
\end{enumerate}
\end{theorem}
\begin{proof}
For parameters $\theta=(y,r_0,\kappa,\mathbf f)$ as in
Proposition~\ref{prop:bounded-type-representation}, let $p_\theta$
denote the function in \eqref{eq:bounded-type-representation} on $B_M$.
Then
\[
 \sup_{x\in B_M}|p_\theta(x)-p_{\tilde\theta}(x)|
 \le (1+\sqrt2)\h(y-\tilde y)+|r_0-\tilde r_0|+|\kappa-\tilde\kappa|
      +\sum_{k=1}^4\|f_k-\tilde f_k\|_\infty .
\]
Indeed, $\h$ is $1$-Lipschitz and $Q$ is $\sqrt2$-Lipschitz in the
Hilbert norm. At a fixed argument, changing $f_k$ changes the maximum
in \eqref{eq:limit-seminorm} by at most $\|f_k-\tilde f_k\|_\infty$.
The reverse triangle inequality in $\mathbb R^3$ gives the estimate,
with the last two terms of the sum omitted if $I_\infty$ is empty.

Suppose $I$ is scattered. A countable order is scattered exactly when
it has only countably many initial segments
\cite[Section~2.2, following Theorem~6]{BruyereCarton2007}.
The formulation in terms of descending cuts is given in
\cite[Exercise~5.33(1)]{Rosenstein1982} and proved in
\cite[Proposition~2.4]{KachMontalban2011}. A proper initial
segment whose complement has a least element is determined by
that element. Thus in a countable order only countably many
initial segments can fail to give a descending cut.
The suborder $I_\infty$ is also scattered. On either order, rounding
a bounded monotone function to a finite mesh gives a uniform
approximation. Each threshold of the rounded function determines
an initial or a final segment. Thus for a fixed finite mesh there
are only countably many rounded functions. Meshes tending to zero
give separability in the supremum metric for each of the four
families of monotone functions above.

The map $\theta\mapsto p_\theta$ is therefore Lipschitz from a separable
product into the bounded functions on $B_M$ with the uniform metric.
Its image contains $S_{\mathrm{norm}}(M)$, so this space is separable.

Now suppose that $I$ contains a copy $(i_r)_{r\in\mathbb Q\cap(0,1)}$
of the rational order. Lupini uses irrational cuts to obtain continuum
many uniformly separated norm types for an additive Hilbert renorming
whose complete Banach theory is NIP
\cite[Corollary~7.4 and Theorem~9.5]{LupiniNormFormulas2026}.
For each irrational $t\in(0,1)$ choose rationals $r_n\uparrow t$ and put
$y_n=E_{i_{r_n}}/\sqrt2$. These vectors and the test vectors
$x_r=F_{i_r}/\sqrt2$ belong to $B_M$ with norm one.
The bounded $1$-Lipschitz functions $p_{y_n}$ have a pointwise
convergent subsequence by diagonal extraction on a countable dense
subset of $B_M$. Write $p_t$ for its limit. The squared Hilbert
norm of $(E_i+F_j)/\sqrt2$ is one. Its seminorm is $\sqrt2$ for
$i<j$ and $1/\sqrt2$ for $i\ge j$, since only in the former case
can the two coefficients occur in the same sum. Thus if $t<r<t'$
with $r$ rational and $t,t'$ irrational, eventual order comparison
gives $p_t(x_r)=\sqrt3$ and $p_{t'}(x_r)=\sqrt{3/2}$. These types
give the separation bound in (\ref{item:strong-types-separation}).

Finally, for $I\models\operatorname{Th}(\mathbb N,<)$ a copy of
$\mathbb Q$ in $I$ misses the initial $\omega$-chain and meets each
remaining finite-distance block at most once. Indeed each point in
that copy has infinitely many predecessors and every two distinct
points have infinitely many points of the copy between them. Its
projection therefore embeds $\mathbb Q$ in $B_I$. Conversely, one
representative from each block of a copy of $\mathbb Q$ in $B_I$
gives a copy in $I$. This proves the criterion for the models of
Theorem~\ref{thm:separable-classification}.
\end{proof}

For $Y=Y_{\mathbb N}^\circ\oplus_2\mathbb R$ the block order is
empty, so its bounded norm types are strongly separable. Taking
$I=\omega+\sum_{r\in\mathbb Q}\mathbb Z$ and $K=\mathbb R$ gives
a separable elementary extension of $Y$ by the initial-chain
embedding of Lemma~\ref{lem:successor-inclusion}. Its bounded
norm types contain the separated family of
Theorem~\ref{thm:strong-types}(\ref{item:strong-types-separation}).

\appendix
\section{Failure of quantifier elimination}\label{app:QE}
\begin{proposition}[Failure of quantifier elimination]\label{cor:no-QE}
The complete theory of $\Y$ fails quantifier elimination both in the
pure Banach language and in its expansion by the predicates $\h$ and
$q$ on the Banach ball sorts.
\end{proposition}
\begin{proof}
Put $a_1=E_1/\sqrt2$ and $a_2=E_2/\sqrt2$. For $i=1,2$, we have
$\Nrm(a_i)=1$ and $\h(a_i)=q(a_i)=1/\sqrt2$. Every parameter-free
term in one vector variable is a scalar multiple of that vector.
Thus the values of $\Nrm$, $\h$, and $q$ on terms and the Banach
distances between terms agree at $a_1$ and $a_2$. The two vectors
have the same quantifier-free type in the expanded language.

The Hilbert norm is definable in the pure Banach language by
Lemma~\ref{lem:intrinsic-h}. Polarization defines the inner product
by $\langle x,e\rangle_{\h}=(\h(x+e)^2-\h(x-e)^2)/4$.
Corollary~\ref{cor:first-row-definable} gives definable quantification
over the signed first row $\widetilde P$. Hence the predicate
\[
 x\longmapsto\sup_{e\in\widetilde P}|\langle x,e\rangle_{\h}|
\]
is definable in the pure Banach language. Every $e\in\widetilde P$
in $\Y$ is of the form $\pm(E_1+F_j)/\sqrt2$ for some $j>1$.
Hence the supremum is $1/2$ at $a_1$ and zero at $a_2$. Their
complete types therefore differ already in the pure Banach language.
\end{proof}

\section{Strict order property}\label{app:SOP}
We use the continuous-logic SOP convention of
\cite[Definition~3.4]{Khanaki2026}. A formula $\varphi(x,y)$
witnesses SOP if there are $\varepsilon>0$ and sequences
$(a_i,b_i)_{i<\omega}$ in a monster model $\mathfrak U$ for which
\begin{equation}\label{eq:Khanaki-SOP}
 \varphi(\mathfrak U,a_i)\le
 \varphi(\mathfrak U,a_{i+1}),\qquad
 \varphi(b_j,a_i)+\varepsilon
 \le\varphi(b_{i+1},a_j)\quad(i<j),
\end{equation}
where the first inequality holds pointwise on the
entire appropriate sort of the monster model.
\begin{theorem}[A quantifier-free SOP witness]\label{thm:SOP}
The complete theory $\operatorname{Th}(\Y)$ has SOP,
witnessed by a quantifier-free formula.
\end{theorem}
\begin{proof}
Let $E_i=(e_i,0)$, $F_j=(0,e_j)$, and let $B_2$ be
the ball sort $\{z:\Nrm(z)\le2\}$.
For a positive rational number $t$ put
\[
 G_t(x,y)=\Nrm(x+ty)^2+\Nrm(x-ty)^2
\]
For $x\in B_2$ and $y=F_j$ with $j\ge3$, any $t>4$ makes the
perturbed column dominate the other columns in both seminorms
$q(x\pm tF_j)$. We use the consecutive integers $5$ and $6$.
In $G_6(x,F_j)-G_5(x,F_j)$ the contribution of the quadratic
term in $t$ is $4(6^2-5^2)=44$. Subtracting this constant cancels
that contribution. Dividing by $1000$ and adding $1/2$ gives
a formula with values in $[0,1]$ on $B_2^2$.
For $x,y\in B_2$ put
\begin{equation}\label{eq:explicit-SOP-formula}
 \varphi(x,y)=
 \frac12+\frac{G_6(x,y)-G_5(x,y)-44}{1000}.
\end{equation}
This formula is quantifier-free because it uses rational scalar
multiplication, addition, norm predicates, and
continuous connectives. For $x,y\in B_2$, the triangle inequality
gives $\Nrm(x\pm ty)\le 2(1+t)$, so the vectors $x\pm ty$ lie in
$B_{12}$ for $t=5$ and in $B_{14}$ for $t=6$.
Thus $0\le G_6\le 2\cdot14^2=392$ and
$0\le G_5\le 2\cdot12^2=288$ on $B_2^2$,
which gives $0.168\le\varphi\le0.848$.
Thus \eqref{eq:explicit-SOP-formula} is a
$[0,1]$-valued formula.

For $j\ge3$ and $x=(u,v)\in B_2$, write
\[
 M_j(x)=\max_{k<j}u_k,\qquad
 m_j(x)=\min_{k<j}u_k.
\]
Since $q(x)\le\Nrm(x)\le2$, all terms
$|u_k+v_j\pm t|$ with $k<j$ are at least
$t-2>2$ for $t\in\{5,6\}$, while terms of
$q(x\pm tF_j)$ with second index different
from $j$ are at most $2$. For each $k<j$ the sum
$u_k+v_j+t$ is positive and $u_k+v_j-t$ is negative.
Taking the maximum of their absolute values over $k<j$
for each sign gives
\[
 q(x+tF_j)=t+v_j+M_j(x),\qquad
 q(x-tF_j)=t-v_j-m_j(x).
\]
For $t\in\{5,6\}$, applying the parallelogram identity to the
Hilbert terms and squaring the two displayed values of $q$ gives
\[
 \begin{split}
 G_t(x,F_j)
 &=2\h(x)^2+4t^2
   +2t(M_j(x)-m_j(x))\\
 &\quad +(v_j+M_j(x))^2+(v_j+m_j(x))^2.
 \end{split}
\]
Subtracting the identities for $t=6$ and $t=5$
gives $G_6-G_5=44+2(M_j-m_j)$. Substituting into
\eqref{eq:explicit-SOP-formula}, we obtain
\begin{equation}\label{eq:explicit-SOP-identity}
 \varphi(x,F_j)=
 \frac12+\frac{M_j(x)-m_j(x)}{500}.
\end{equation}

For $i<\omega$ choose
\[
 a_i=F_{i+3},\qquad b_i=E_{i+2}.
\]
These vectors lie in $B_2$, since
$\Nrm(E_k)=\Nrm(F_j)=\sqrt2$ for $k\ge1$,
$j\ge2$.  The quantities $M_j(x)$ are
nondecreasing and $m_j(x)$ nonincreasing
as $j$ increases.  Therefore
\eqref{eq:explicit-SOP-identity} shows that
\[
 \varphi(x,a_i)\le\varphi(x,a_{i+1})
 \qquad\text{for every }x\in B_2.
\]

For $i<j$, all first-component coordinates of $b_j=E_{j+2}$
with $k<i+3$ are zero. Hence $M_{i+3}(b_j)=m_{i+3}(b_j)=0$
and $\varphi(b_j,a_i)=1/2$.
For $b_{i+1}=E_{i+3}$, the coordinates with $k<j+3$ include
a $1$ at $k=i+3$ and a $0$ at $k=1$. Thus
$M_{j+3}(b_{i+1})=1$ and $m_{j+3}(b_{i+1})=0$, giving
$\varphi(b_{i+1},a_j)=1/2+1/500$.
Thus
\[
 \varphi(b_j,a_i)+\frac1{500}
 =\varphi(b_{i+1},a_j)\qquad(i<j).
\]
For each fixed $i$, the formula
\[
 \sup_{x\in B_2}
  (\varphi(x,a_i)-\varphi(x,a_{i+1}))_+
\]
has value zero in $\Y$. In a monster elementary extension it
still has value zero, so the pointwise inequality holds on the
entire ball $B_2$ of that extension. The inequalities on
$a_i,b_j,b_{i+1},a_j$ retain the margin $1/500$. Thus
$\varphi$ witnesses SOP in $\operatorname{Th}(\Y)$.
\end{proof}

\section*{Declaration on the use of AI}
The authors made extensive use of ChatGPT (OpenAI), with the
GPT-6.1~Sol, GPT-6.2~Sol, and GPT-6~Astra models, while developing
and revising this manuscript. In response to a question about
approximating formula values, ChatGPT suggested an approach that helped
shape the proof of NIP for the complete theory of $Y$. It also suggested
mathematical ideas and proof strategies, helped develop parts of several
arguments, and was especially helpful with
Lemma~\ref{lem:relative-positions} and
Theorem~\ref{thm:separable-classification}. The authors used the models
to look for gaps and clarify the relationship between their NIP proof
and Lupini's work. The authors directed the investigation, evaluated
and selected the arguments, and checked and revised the proofs. They
take full responsibility for the mathematical statements, proofs, and
exposition.

\bibliographystyle{abbrvurl}
\bibliography{bibliography}
\end{document}